\documentclass[11pt]{article}

\def\mathbb{\Bbb}
\usepackage{amsfonts,latexsym,amstext}
\usepackage{amsmath}
\usepackage{amssymb}

\usepackage{comment}
\usepackage{amsthm}
\usepackage[latin1]{inputenc}

\theoremstyle{plain}
\newtheorem{theorem}{Theorem}[section]
\newtheorem{lemma}[theorem]{Lemma}
\newtheorem{proposition}[theorem]{Proposition}
\newtheorem{corollary}[theorem]{Corollary}
\newtheorem{definition}{Definition}[section]

\theoremstyle{hypothesis}
\theoremstyle{remark}
\newtheorem{remark}[theorem]{Remark}

\makeatletter
\@addtoreset{equation}{section}

\makeatother

\def\qed{{\hfill\hbox{\enspace${ \square}$}} \smallskip}
\def\sqr#1#2{{\vcenter{\vbox{\hrule height .#2pt \hbox{\vrule
 width .#2pt height#1pt \kern#1pt \vrule
width .#2pt} \hrule height .#2pt}}}}
\def\square{\mathchoice\sqr54\sqr54\sqr{4.1}3\sqr{3.5}3}

\def\ds{\begin{displaystyle}}
\def\eds{\end{displaystyle}}
\def\dis{\displaystyle }
\def\<{\langle }
\def\>{\rangle }

\def\esssup{\mathop{\rm esssup}}

\def\R{\mathbb R}
\def\N{\mathbb N}

\def\E{\mathbb E}
\def\P{\mathbb P}

\def\F{\mathbb F}

\def\cala{{\cal A}}
\def\calb{{\cal B}}

\def\cale{{\cal E}}
\def\calf{{\cal F}}
\def\calg{{\cal G}}

\def\calp{{\cal P}}

\def\call{{\cal L}}

\newcommand{\spernu}[1]{\mathbb{E}_{\nu}^{t,x,a} \left[ #1 \right]}
\newcommand{\spertxanuepsilon}[1]{\mathbb{E}^{t,x,a}_{\nu^{\epsilon}} \left[ #1 \right]}                               % per la speranza sotto P_nu_ast
\newcommand{\spertxa}[1]{\mathbb{E}^{t,x,a} \left[ #1 \right]}                               % per la speranza
\DeclareMathAlphabet{\mathonebb}{U}{bbold}{m}{n}                           %
\newcommand{\one}{\ensuremath{\mathonebb{1}}}                               % per le funzioni indicatrici
\title{
Constrained BSDEs representation of the value function
in optimal control of pure jump Markov processes}
\date{}
\author{Elena Bandini and Marco Fuhrman
\\Politecnico di Milano,
Dipartimento di Matematica\\
via Bonardi 9, 20133 Milano, Italy\\
e-mail: elena.bandini@polimi.it, marco.fuhrman@polimi.it}

\begin{document}

\maketitle

\begin{abstract}
We consider a classical finite horizon optimal control problem
for  continuous-time pure jump Markov processes described by
means of a rate transition measure depending on a control parameter
and controlled by a feedback law. For this class of problems
the value function can often be described as  the unique solution to
the corresponding Hamilton-Jacobi-Bellman equation.
We prove a probabilistic representation
for the   value function, known as
nonlinear Feynman-Kac formula. It relates the
value function with a backward stochastic differential equation
 (BSDE) driven by a random measure and with a sign constraint on
 its martingale part.
We also prove existence and uniqueness results  for this class of constrained BSDEs.
  The connection of the control problem with the constrained BSDE uses a
 control randomization method
 recently de\-ve\-loped
 in the works of I. Kharroubi and H. Pham and their co-authors.
This  approach also allows to prove that the value function of the original non-dominated control problem coincides with the value function of an auxiliary dominated control problem, expressed in terms of equivalent changes of probability measures.
\end{abstract}

\section{Introduction}

The main aim of this paper is to prove that the value function
in a classical optimal control problem for pure jump Markov processes
can be represented by means of an appropriate Backward Stochastic Differential
Equation (BSDE) that we introduce and
for which we prove an existence and uniqueness result.

We start by   describing our setting in an informal way. A
 pure jump Markov process $X$ in a general measurable
state space $(E,\cale)$ can be described by means of
a rate transition measure, or intensity measure,
  $ \nu(t,x,B)$ defined for $t\ge 0$, $x\in E$, $B\in \cale$.
  The process starts at
  time $t\ge 0$ from some initial point  $x\in E$
  and stays there up to a random time $T_1$
  such that
  $$
  \P (T_1>s)=\exp\left( -\int_t^s \nu(r,x,E)\,dr\right),
  \qquad s\ge t.
  $$
At time $T_1$, the process jumps to a new point $X_{T_1}  $
chosen with probability $\nu(T_1,x,\cdot)/\nu(T_1,x,E)$ (conditionally to $T_1$)
and then it stays again at $X_{T_1}  $ up to another random
 time $T_2$ such that
 $$
  \P (T_2>s\mid T_1, X_{T_1})=\exp\left( -\int_{T_1}^s \nu(r,X_{T_1},E)\,dr\right),
  \qquad s\ge {T_1},
  $$
and so on.

A controlled pure
jump Markov process is obtained starting from a
 rate measure
  $ \lambda(x,a,B)$  defined for  $x\in E$, $a\in A$, $B\in \cale$,
  i.e.,  depending on a control parameter $a$ taking values in a
  measurable space of control actions $(A,\cala)$.
A natural way to control   a Markov process is to choose a feedback
control law, which is a measurable function $\alpha: [0,\infty)\times E\to A$.
$\alpha(t,x)\in A$ is the control action selected at time $t$ if the system
is in state $x$.  The controlled Markov process $X$ is simply
the one corresponding to the rate transition measure
  $ \lambda(x,\alpha(t,x),B)$. Let us denote by $\P^{t,x}_\alpha$
  the corresponding law, where $t,x$ are the initial time and starting point.

We note that an alternative construction of  (controlled
or uncontrolled) Markov processes consists in defining them as
 solutions to stochastic equations driven by
some noise (for instance, by a Poisson process) and with appropriate  coefficients
depending on a control process.
  In the context of pure jump processes, our approach based on the
  introduction of the controlled
rate measure
  $ \lambda(x,a,B)$ often leads to
more general results and it is more natural in several contexts.

In the classical finite horizon control problem one seeks to maximize
over all control laws  $\alpha$
a functional of the form
\begin{equation}\label{functional_cost_intro}
J(t,x,\alpha)= \mathbb{E}^{t,x}_{\alpha}\left[ \int_{t}^{T} f(s,\,X_{s}, \alpha(s,X_s))\, ds
+g(X_{T}) \right],
\end{equation}
where a deterministic finite horizon $T>0$ is given and $f,g$ are given real functions,
defined on $ [0,T]\times E\times A$ and $ E$, representing the running
cost and the terminal cost, respectively.
The value function of the control problem is defined in the usual way:
\begin{equation}\label{value_function_intro}
v(t,x) = \sup_{ \alpha } J(t,x,\alpha), \qquad t\in [0,T], \,x\in E.
\end{equation}

We will only consider the case when the controlled rate measure $\lambda$
and the costs $f,g$ are bounded. Then, under some technical assumptions,
$v$ is known to be the unique solution on $[0,T]\times E$
to the Hamilton-Jacobi-Bellman (HJB)
equation
\begin{equation}\label{HJB_intro}
\left\{\begin{array}{lll}\dis
-\frac{\partial v}{\partial t}(t,x)&=&\dis \sup_{a \in A}\left(
 \int_{E}(v(t,y)-v(t,x))\,\lambda(x,a,dy)  + f(t,x,a)\right) ,
\\
v(T,x)&=& g(x),
\end{array}\right.
\end{equation}
and if the supremum is attained at some $\alpha(t,x)\in A$ depending
measurably on $(t,x)$ then $\alpha$ is an optimal feedback law.
Note that the right-hand side of \eqref{HJB_intro} is an integral
operator: this allows for easy notions of solutions to the HJB equation,
that do not in particular need the use of the theory of viscosity solutions.

Our purpose is to relate the value function $v(t,x)$ to
 an appropriate BSDE.  We wish to extend to our
framework the theory developed in the context of classical optimal control for
diffusion processes, constructed as solutions to stochastic differential equations
of Ito type driven by Browian motion, where   representation formulae
for  the solution to the HJB equation
exist  and are often called non-linear Feyman-Kac formulae.
The  majority of those results requires that only the drift coefficient
of the stochastic equation depends on the control parameter,
so  that in this case the HJB equation is a second-order semi-linear
partial differential equation and the non-linear Feyman-Kac formula
is well known, see e.g. \cite{EPQ}. Generally, in this case
the laws of the corresponding controlled processes are all absolutely continuous
with respect to the law of a given, uncontrolled process, so that they form
a dominated model.

A natural extension to our framework could be obtained imposing conditions implying that
the set of probability laws $\{\P^{t,x}_\alpha\}_\alpha$, when $\alpha$
varies over all feedback laws, is a dominated model. This is the point of view taken
in \cite{CoFu-m}, where an appropriate BSDE is
introduced and solved and a Feyman-Kac formula for the value function is proved
in a restricted framework. Extensions are given in
\cite{BaCo}
to controlled semi-Markov processes and in
\cite{CoFu-mpp}
to more general non-Markovian cases.

In the present paper we want to consider the general case when
$\{\P^{t,x}_\alpha\}_\alpha$ is not a dominated model. Even
for finite state space $E$, by a proper
choice of the measure $\lambda(x,a,B)$ it is easy to
formulate quite natural control problems for which this is the case.

In the context of controlled
diffusions,  probabilistic formulae for the value function
 for non-dominated models have been discovered only in recent years. We note
 that in this case the HJB equation is a fully non-linear partial differential equation.
 To our knowledge, there are only a few available techniques.
One possibility is to
use the theory of second-order BSDEs, see for instance \cite{ChSoToVi07},
\cite{SoToZh11}.
Another possibility relies on the use of the theory of $G$-expectations,
see e.g. \cite{Pe06}.  Both theories have been largely developed by several authors.
In this paper we rather follow another approach which is presented in the
paper \cite{KhPh} and was predated by similar results concerning
optimal switching  or optimal impuse control problems,
see \cite{ElKh10}, \cite{ElKh14},
\cite{ElKh14a}, \cite{KhMaPhZh}, and followed by some
extensions and applications, see  \cite{FuPh}, \cite{Co-Chou},
\cite{CoFuPh}.
It consists in a {\em control randomization
method} (not to be confused with the use of relaxed controls)
which can be described informally as follows, in our framework of
controlled pure jump Markov processes.

We note that for any choice of a  feedback law $\alpha$  the
 pair of stochastic processes $(X_s, \alpha(s,X_s))$
 represents the state trajectory and the associated control process. In a first step,
 for any initial time $t\ge 0$ and starting point $x\in E$, we replace
 it  by an (uncontrolled) Markovian pair of pure jump
 stochastic processes $(X_s, I_s)$, possibly
 constructed on a different probability space, in such a way that
 the process $I$ is a Poisson process with values in the space of
 control actions $A$ with an  intensity measure $\lambda_0(da)$
which is arbitrary but finite and with full support.
Next we formulate an auxiliary
optimal control problem where we control
the intensity of the process $I$: for any predictable, bounded and positive  random field
$\nu_t(a)$, by means of a theorem of Girsanov type
we construct a probability measure $\P_\nu$ under
which the compensator
of $I$ is the random measure
$\nu_t(a)\, \lambda_0(da)\,dt$ (under $\P_\nu$ the law of $X$ also changes)
and then we maximize the functional
$$
 \mathbb{E}_{\nu}\left[g(X_{T}) + \int_{t}^{T} f(s,\,X_{s},\,I_{s}) \,ds\right],
$$
over all possible choices of the process $\nu$. Following the
terminology of \cite{KhPh}, this will be called the {\em dual}
control problem. Its value function, denoted
    $v^*(t,x,a)$, also depends  {\em a priori} on the starting point $a\in A$
    of the process $I$  (in fact we should  write
    $\P_\nu^{t,x,a}$ instead of $\P_\nu$, but in this discussion
    we drop this dependence for simplicity)
    and the family $\{\P_\nu\}_\nu$ is a dominated model.
As in \cite{KhPh} we are able
    to show that the value functions for the original problem and the dual one
are the same: $v(t,x)=v^*(t,x,a)$, so that the latter does not in fact
depend on $a$. In particular we have replaced the original control problem
by a dual one
that corresponds to a dominated model and has the same value function.
Moreover, we can introduce a well-posed BSDE that represents
$v^*(t,x,a)$ (and hence $v(t,x)$). It is an equation
on the time interval $[t,T]$ of the form
\begin{eqnarray}\label{BSDE_intro}
Y_{s} &=& g(X_T) + \int_{s}^{T}f(r,X_r,I_r)\,dr + K_T - K_s
\nonumber\\
&& - \int_{s}^{T}\int_{E \times A}Z_r(y,\,b)\, q(dr\,dy\, db)
- \int_{s}^{T}\int_{A}Z_r(X_r,\,b)\, \lambda_0(db)\,dr,
\end{eqnarray}
with unknown triple $(Y,Z,K)$ (depending also on $(t,x,a)$),
where $q$ is the compensated random measure
associated to $(X,I)$, $Z$ is a predictable random field and $K$ a predictable increasing
càdlàg process, where
we additionally add the sign constraint
\begin{equation}\label{BSDE_constraint_intro}
Z_s(X_{s-},b)\leqslant 0.
\end{equation}
It turns out that this equation has a unique minimal solution, in an appropriate
sense, and that the value of the process $Y$ at the initial time
represents both the original and the dual value function:
\begin{equation}\label{FK_intro}
    Y_t= v(t,x)=v^*(t,x,a).
\end{equation}
This is the desired BSDE representation of the value function for the
original control problem and a Feyman-Kac formula for the general HJB equation
\eqref{HJB_intro}.

The paper is organized as follows. Section \ref{Sec_preliminaries} is essentially
devoted to lay down a setting where the classical optimal control
problem \eqref{value_function_intro} is solved by means
of  the corresponding  HJB equation \eqref{HJB_intro}.
We first recall the general construction of a Markov process given its rate transition measure.
Having in mind to apply techniques based on BSDEs driven by random
measures we need to work in a canonical
setting and use a specific filtration, see Remark
\ref{filtrazionidiverse}. Therefore the construction we present is based on the
well-posedness of the martingale problem for  multivariate (marked) point processes
studied in \cite{J_a} and it is exposed in detail.
This general construction is then used
to formulate in a precise way the optimal control problem for
the jump Markov process and it is
used again in the subsequent section when we define
the pair $(X,I)$ mentioned above.
Still in section \ref{Sec_preliminaries}, we present classical results
on existence and uniqueness
of the solution to the HJB equation \eqref{HJB_intro} and its
identification with the value function $v$. These results are similar to those in
\cite{Pl}, a place where we could find a clear and complete exposition of all
the basic theory and to which we refer for further references
and related results.
We note that the compactness of the space of control actions
$A$, together with suitable upper-semicontinuity conditions of the
coefficients of the control problem, is one of the standard assumptions needed
to ensure the existence of an optimal control, which is usually constructed
by means of an appropriate measurable selection theorem.
Since our main aim was only to find
a representation formula for the value function we wished to avoid
the compactness condition. This was made possible by the use of a different
measurable selection result, that however requires lower-semicontinuity conditions.
Although this is not usual in the context of maximization problems,
this turned out to be the right condition that allows to dispense
with compactness assumptions and to prove well-posedness of the HJB equation
and a verification theorem. A small variation of the proofs recovers
the classical results in \cite{Pl}, and even with slightly weaker assumptions:
see Remark \ref{confrontoPliska} for a more detailed comparison.

In section \ref{Section_control_rand} we start to develop the control randomization
method: we introduce the auxiliary process $(X,I)$ and formulate the dual
control problem under appropriate conditions. Finding the correct formulation
required some efforts; in particular we could not mimic the approach of previous
works on control randomization mentioned above,
since we are not dealing
with processes defined as solutions to stochastic equations.

In section \ref{Section_BSDE_rep}
we introduce the constrained BSDE
\eqref{BSDE_intro}-\eqref{BSDE_constraint_intro} and we prove, under suitable
conditions, that it has a unique minimal solution
$(Y,Z,K)$ in a certain class of processes. Moreover, the value of $Y$
at the initial time coincides with the  value function of the dual optimal
control problem.
This is the content of the first
of our main results, Theorem
\ref{Thm_ex_uniq_minimal_BSDE}. The proof relies on a penalization
approach and a monotonic passage to the limit,
and combines BSDE techniques with control-theoretic arguments:
for instance, a ``penalized'' dual control problem is also introduced in
order to obtain certain uniform upper bounds. In \cite{KhPh},
in the context of diffusion processes, a more general result is proved,
in the sense that the generator $f$ may also depend on $(Y,Z)$;
similar generalizations are possible in our context as well, but they seem
less motivated and in any case they are
not needed for the applications to optimal control.

Finally, in section \ref{FeynmanKac} we prove the second
of our main results, Theorem  \ref{thm_conv_vn}. It states that the initial value of
the process $Y$ in \eqref{BSDE_intro}-\eqref{BSDE_constraint_intro}
 coincides with the value function $v(t,x)$.
As a consequence, the value function is the same for the original optimal
control problem and for the dual one and we have the non-linear
Feynman-Kac formula \eqref{FK_intro}.

The assumptions in Theorem  \ref{thm_conv_vn}
are fairly general: the state space $E$ and the control action space $A$
are Borel spaces, the controlled kernel $\lambda$ is bounded
and has the Feller property, and the cost functions $f,g$ are continuous and bounded.
No compactness assumption is required.
When $E$ is finite or countable we have the special case of (continuous-time) controlled
Markov   chains. A large class of optimization problems for controlled Markovian
queues falls under the scope of our result.

In recent years there has been much interest in numerical approximation
of the value function in optimal control of
Markov processes, see for instance the book \cite{GuHe} in the discrete
state case.
The Feynman-Kac formula \eqref{FK_intro} can be used to design
algorithms based on numerical approximation of the solution to
the constrained BSDE \eqref{BSDE_intro}-\eqref{BSDE_constraint_intro}.
Numerical schemes for this kind of equations have been proposed
and analyzed in the context of diffusion processes, see
\cite{KhLaPh}, \cite{KhLaPha}.  We hope that our results may
be used as a foundation for similar methods in the context of pure jump
processes as well.

\section{Pure jump controlled Markov processes}\label{Sec_preliminaries}

\subsection{The construction of a jump Markov process given the rate transition measure}\label{Section_construction_X}

Let $E$ be a Borel space, i.e., a topological space homeomorphic
to a Borel subset of a compact metric space (some authors call it a Lusin space);
in particular, $E$ could be a Polish space.
Let $\mathcal{E}$ denote the corresponding Borel $\sigma$-algebra.

We will often need to construct a Markov process in $E$ with a given
(time dependent)
rate transition measure, or intensity measure, denoted by
$\nu$.
With this terminology we
mean that $B\mapsto \nu(t,x,B)$ is a nonnegative measure on $(E, \mathcal{E})$ for every
$(t,x)\in [0,\infty)\times E$ and $(t,x)\mapsto \nu(t,x, B)$ is a Borel measurable function
on $[0,\infty) \times E$ for every $B \in  \mathcal{E}$.
We assume that
\begin{equation}\label{nulimitato}
\sup_{t\ge0,\,x \in E}\, \nu(t,x,E) < \infty.
\end{equation}
We recall the main steps in the construction of the corresponding Markov
process. We note that \eqref{nulimitato} allows to construct a non-explosive
process. Since $\nu$ depends on time the process will not be
time-homogeneous in general. Although the existence of
such a process is a well known fact, we need special care in the choice
of the corresponding filtration, since this will be crucial when we solve
associated BSDEs and implicitly apply
a version of the martingale representation theorem in the sections that follow:
see also Remark \ref{filtrazionidiverse} below. So in the following we
will use an explicit construction that we are going to describe. Many of the
techniques we are going to use are borrowed from the  theory of
multivariate (marked) point processes. We will often follow
\cite{J_a}, but we also refer the reader to the treatise  \cite{BrLa}
for a more systematic exposition.

We start by constructing a suitable sample space to describe
the jumping mechanism of the Markov process.
 Let $\Omega'$ denote the set of sequences $\omega'=(t_n,e_n)_{n \geq 1}$
in $ ((0,\infty)\times E )\cup \{(\infty,\Delta)\}$, where $\Delta\notin E$
 is   adjoined to  $E$ as an isolated point, satisfying in addition
 \begin{equation}\label{processo_punto}
    t_n \le t_{n+1}; \qquad
  t_n< \infty\;\Longrightarrow\; t_n < t_{n+1}.
 \end{equation}
To describe  the initial condition we will use the measurable space
$(E,\cale)$.
Finally, the sample space for the Markov process will be
$\Omega=E\times \Omega'$. We define canonical
functions $T_n:\Omega\to (0,\infty]$,
$E_n:\Omega\to E\cup\{\Delta\}$
as follows:  writing $\omega=(e,\omega')$ in the form
$\omega=(e,t_1,e_1,t_2,e_2,\ldots)$ we set  for $t\ge 0$ and for $n\ge 1$
$$
T_n (\omega )= t_n ,
\qquad E_n(\omega)=e_n,
\qquad T_\infty  (\omega )= \lim_{n\to\infty} t_n,
\qquad T_0(\omega)=0, \qquad E_0(\omega)=e.
 $$
 We also define
 $X:\Omega\times [0,\infty)
\to E\cup\{\Delta\}$ setting
  $X_t = 1_{[0, T_1]}(t)\,E_0+\sum_{n\ge 1}
  1_{ (T_n,T_{n+1}]  }(t)\,E_n$ for $t<T_\infty$, $X_t=\Delta$ for $t\ge T_\infty$.

 In $\Omega$ we introduce for all $t\ge0$ the $\sigma$-algebras
 $\calg_t=\sigma(N(s,A)\,:\, s\in (0,t],A\in\cale)$, i.e. generated
 by the counting processes defined as $N(s,A)=
 \sum_{n\ge 1}1_{T_n\le s}1_{E_n\in A}$.

To take into account the initial condition we also introduce
the filtration $\F=(\calf_t)_{t\ge 0}$, where
$\calf_0=\cale\otimes \{\emptyset,\Omega'\}$,
and for all $t\ge0$
$\calf_t$ is the $\sigma$-algebra generated by $\calf_0$ and $\calg_t$.
$\F$ is right-continuous and will be called the natural filtration.
In the following all concepts of measurability for stochastic processes
(adaptedness, predictability etc.) refer to $\F$.
 We denote by $\calf_\infty$
the $\sigma$-algebra generated by all the $\sigma$-algebras $\calf_t$.
The symbol $\calp$ denotes the $\sigma$-algebra
of $\mathbb{F}$-predictable subsets of $[0,\infty) \times \Omega$.

The initial distribution of the process $X$ will be described by
a probability measure $\mu$ on $(E,\cale)$.
Since
$\calf_0=\{A\times \Omega'\,:\,A\in\cale\}$ is isomorphic
to $\cale$,
$\mu$ will be identified with a probability measure on $\calf_0$, denoted by the same
symbol (by abuse of notation) and such that $\mu(A\times \Omega')=\mu(A)$.

On the filtered sample space
$(\Omega,\F)$ we have so far introduced
the canonical marked point process $(T_n,E_n)_{n\geq1}$.
 The corresponding
 random measure $p$ is, for any $\omega\in\Omega$,
 a $\sigma$-finite measure on $((0,\infty)\times E,\calb((0,\infty))\otimes \cale)$
  defined as
\begin{eqnarray*}
p(\omega, ds\, dy)= \sum_{n\ge 1}1_{ T_n(\omega) < \infty}\,\delta_{(T_n(\omega),E_n(\omega))}(ds\,dy),
\end{eqnarray*}
where $\delta_{k}$ denotes the Dirac measure at point $k\in (0,\infty)\times E$.

Now let $\nu$ denote a time-dependent
rate transition measure
 as before, satisfying \eqref{nulimitato}.  We need
 to introduce the corresponding generator and transition semigroup as follows.
 We denote by $B_b(E)$ the space
 of $\cale$-measurable bounded real functions on $E$ and for
 $\phi\in B_b(E)$ we set
\begin{displaymath}
\mathcal{L}_t\phi(x) = \int_{E}(\phi(y)-\phi(x))\,\nu(t,x,dy), \qquad t\ge0,\,x \in E.
\end{displaymath}
For any $T \in (0,\infty)$ and  $g\in B_b(E)$ we consider the Kolmogorov equation
on $[0, T]\times E$:
\begin{equation}\label{kolm_eq}
\left\{\begin{array}{l}
\dis \frac{\partial v}{\partial s}(s,x) + \mathcal{L}_sv(s,x)  =  0 ,
 \\
v(T,x) = g(x).
\end{array}
\right.
\end{equation}
It is easily proved
that there exists a unique measurable bounded function $v: [0,\,T] \times E$ such that
$v(T,\cdot)=g$ on $E$ and,
  for all $x \in E$, $s \mapsto v(s,x)$ is an absolutely continuous map on $[0,T]$ and
  the first equation in \eqref{kolm_eq} holds for almost all $s\in [0,T]$ with respect
  to the Lebesgue measure.
 To verify this we first write  \eqref{kolm_eq} in the equivalent
integral form
$$
v(s,x) = g(x) + \int_s^T\mathcal{L}_rv(r,x)\,dr,\qquad s \in [0,T],\; x \in E.
$$
Then, noting the inequality
$|\call_t\phi(x)|\le 2 \sup_{y \in E}|\phi(y)|\, \sup_{t\in[0,T],y \in E}\, \nu(t,y,E) $,
a solution to the latter equation can be obtained by a standard
fixed point argument in the space of
bounded measurable real functions on $[0,\,T] \times E$
endowed with the supremum norm.

 This allows to define the transition operator $P_{sT}:B_b(E)\to B_b(E)$,
for  $0\le s \le T$, letting $P_{sT}[g](x)=v(s,x)$, where $v$ is the
solution to \eqref{kolm_eq} with terminal condition $g\in \calb_b(E)$.

\begin{proposition}\label{prop_construction_X_Markov_pure_jump_1}
Let   \eqref{nulimitato} hold and let us fix $t\in [0,\infty)$
and a probability measure $\mu$ on $(E,\cale)$.
\begin{enumerate}
  \item
There exists a unique probability measure on $(\Omega, \mathcal{F}_\infty)$, denoted by $\mathbb{P}^{t,\mu}$, such that its restriction to $\calf_0$ is $\mu$
and the $\F$-compensator (or dual predictable projection) of the measure
$p$ under $\mathbb{P}^{t,\mu}$ is the random measure
$\tilde{p}(ds\,dy) :=1_{[t,T_\infty)}(s) \nu(s,X_{s-},dy)\,ds$.
Moreover, $\mathbb{P}^{t,\mu}(T_\infty=\infty)=1$.

  \item
  In the probability space $\{\Omega, \mathcal{F}_{\infty},\mathbb{P}^{t,\mu}\}$
  the   process $X$  has distribution $\mu$ at time $t$   and it is
Markov on the time interval $[t,\infty)$ with
respect to $\mathbb{F}$ with transition operator $P_{sT}$:
explicitly,
for every $ t\le s\le T$  and for every $g\in \calb_b(E)$,
\begin{displaymath}
\mathbb{E}^{t,\mu}\left[g(X_T)\mid \mathcal{F}_s\right] = P_{sT}[g](X_s),
\qquad \mathbb{P}^{t,\mu}-a.s.
\end{displaymath}

\end{enumerate}
\end{proposition}

\proof
Point 1 follows from a direct application
of  \cite{J_a}, Theorem 3.6.  The non-explosion
condition $\mathbb{P}^{t,\mu}(T_\infty=\infty)=1$ follows from
the fact that $\lambda$ is bounded.

 To prove point  2 we denote $v(s,x)=P_{sT}[g](x)$ the solution to the
Kolmogorov equation \eqref{kolm_eq} and note that
$$
v(T,X_T)- v(s,X_s) =   \int_{s}^T  \frac{\partial v}{\partial r}(r,X_r)\,dr
+  \int_{(s,T]}\int_{E} (v(r,y) - v(r,X_{r-}))\,p(dr\,dy).
$$
This identity is easily proved taking into account that $X$ is constant among
jump times and using the definition of the random measure $p$.
Recalling the form of the $\mathbb{F}$-compensator $\tilde p$ of $p$ under
$\mathbb{P}^{t,\mu}$ we have, $\mathbb{P}^{t,\mu}$-a.s.,
$$
\begin{array}{l}\dis
\mathbb{E}^{t,\mu}\Big[
\int_{(s,T]}\int_{E} (v(r,y) - v(r,X_{r-}))\,p(dr\,dy)\mid \calf_s\Big]
\\\qquad \dis
= \mathbb{E}^{t,\mu}\Big[
\int_{(s,T]}\int_{E} (v(r,y) - v(r,X_{r-}))\,\tilde p(dr\,dy)\mid \calf_s \Big]
\\\qquad \dis
= \mathbb{E}^{t,\mu}\Big[
\int_{(s,T]}\int_{E} (v(r,y) - v(r,X_{r}))\,\nu(r,X_r,dy)\,dr\mid \calf_s \Big]
\\\qquad \dis
= \mathbb{E}^{t,\mu}\Big[
\int_{(s,T]}\call_r v(r,X_{r}) \,dr\mid \calf_s \Big]
\end{array}
$$
and
we finally obtain
$$
\begin{array}{l}\dis
\mathbb{E}^{t,\mu}\left[g(X_T)\mid \mathcal{F}_s\right] - P_{sT}[g](X_s)
=
\mathbb{E}^{t,\mu}[v(T,X_T)\mid \calf_s ]
- v(s,X_s)
\\\qquad \dis
= \mathbb{E}^{t,\mu}\Big[  \int_{s}^T \Big( \frac{\partial v}{\partial r}(r,X_r)
+\call_r v(r,X_{r})\Big)\,dr
\mid \calf_s \Big] =0.
\end{array}
$$
 \endproof

In the following we will mainly consider initial distributions $\mu$ concentrated
at some point $x\in E$, i.e. $\mu=\delta_x$. In this case we use the notation
$\mathbb{P}^{t,x}$ rather than $\mathbb{P}^{t,\delta_x}$. Note that, $\mathbb{P}^{t,x}$-a.s.,
we have $T_1>t $ and therefore $X_s=x$ for all $s\in [0,t]$.

\begin{remark}
\label{filtrazionidiverse}
Since the process $X$ is $\F$-adapted, its natural filtration
$\F^X=(\calf^X_t)_{t\ge 0}$ defined by
$\calf^X_t =\sigma (X_s\,:\,s\in [0,t])$ is smaller than
$\F$. The inclusion may be strict,  and may remain such
if we consider the corresponding completed filtrations.
The reason is that the random
variables $E_n$ and $E_{n+1}$ introduced above may coincide
on a set of positive  probability, for some $n$, and therefore knowledge of a trajectory of $X$
does not allow to reconstruct the trajectory $(T_n,E_n)$.

In order to have $\calf_s=\calf_s^X$ up to
$\P^{t,\mu}$-null sets one could require that $\nu(t,x,\{x\})=0$,
i.e. that $T_n$ are in fact jump times of $X$, but this would impose unnecessary
restrictions in some constructs that follow.

 Clearly, the Markov property with respect to  $\F$ implies the Markov property
 with respect to  $\F^X$ as well.
\end{remark}

\subsection{Optimal control of pure jump Markov processes}\label{Sect_control_problem}

In this section we formulate and solve an optimal control problem
for a Markov process with a state space $E$, which is still
assumed to be a Borel space with its Borel $\sigma$-algebra $\cale$.
The other data of the problem
will be  another Borel space $A$, endowed with its Borel $\sigma$-algebra
$\cala$ and called the space of control actions;  a finite time horizon, i.e.
a (deterministic) element $T\in (0,\infty)$;
two real valued functions $f$ and $g$, defined on
$[0,T]\times E\times A$ and $E$ and called   running and terminal cost functions
respectively;
and finally a measure transition kernel  $\lambda$  from $(E \times A, \cale\otimes \cala)$
to $(E, \mathcal{E})$:
namely $B\mapsto \lambda(x,a,B)$ is a nonnegative measure on $(E, \cale)$
for every $(x,a)\in E\times A$ and $(x,a)\mapsto \lambda(x,a,B)$ is a Borel measurable function
for every $B\in  \mathcal{E}$.
We assume that $\lambda$ satisfies the following condition:
\begin{equation} \label{lambda_unif_finite}
 \sup_{x \in E, a \in A}\, \lambda(x,a,E) < \infty.
\end{equation}

The requirement that $\lambda(x,a,\{x\})=0$   for all $x \in E$ and
$a\in A$ is natural in many applications, but it is not needed.
The kernel $\lambda$ depending on the control parameter $a \in A$  plays the role of a
controlled  intensity measure for a controlled Markov  process.
Roughly speaking, we may control the dynamics of the process  by changing its jump intensity
dynamically. For a more precise definition, we first
 construct $\Omega$, $\F=(\calf_t)_{t\ge 0}$, $\calf_\infty$
as in the previous paragraph.
Then we introduce the class
of admissible control laws $\mathcal{A}_{ad}$  as the set of
all Borel-measurable maps $\alpha: [0,T]\times E\to A$.
To any such $\alpha$ we associate the rate transition
measure $\nu^\alpha(t,x,dy):=\lambda(x,\alpha(t,x),dy)$.

For every starting time $ t\in [0,T]$ and starting point $x\in E$, and
for each $\alpha\in \mathcal{A}_{ad}$,
we construct as in the previous paragraph the probability measure
on $(\Omega,\calf_\infty)$, that will be denoted
$\mathbb{P}^{t,x}_{\alpha}$,
corresponding to $t$, to the initial distribution concentrated at $x$ and to the
 the rate transition
measure $\nu^\alpha$.
According to
Proposition \ref{prop_construction_X_Markov_pure_jump_1},
 under
$\mathbb{P}^{t,x}_\alpha$ the process $X$ is Markov with respect
to $\F$ and satisfies $X_s = x$ for every $s\in [0,T]$; moreover,
the restriction
 of the measure $p$ to $(t,\infty)\times E$ admits the compensator
$\lambda(X_{s-},\,\alpha(s,X_{s-}),\,dy) \, ds$.
Denoting by $\mathbb{E}^{t,x}_{\alpha}$ the expectation under $\mathbb{P}^{t,x}_{\alpha}$
we finally define,
for $t \in [0,\,T]$, $x \in E$ and $\alpha \in \mathcal{A}_{ad}$,
 the gain functional
\begin{equation}\label{functional_cost}
J(t,x,\alpha)= \mathbb{E}^{t,x}_{\alpha}\left[ \int_{t}^{T} f(s,\,X_{s}, \alpha(s,X_s))\, ds
+g(X_{T}) \right],
\end{equation}
and the value function of the control problem
\begin{equation}\label{value_function}
V(t,x) = \sup_{ \alpha \in \mathcal{A}_{ad}} J(t,x,\alpha).
\end{equation}
Since we will assume below that $f$ and $g$ are at least Borel-measurable
and bounded, both $J$ and $V$ are well defined and bounded.

\begin{remark}
In this formulation the only control strategies
that we consider are control laws of feedback type, i.e.,
the control action $\alpha(t,x)$ at time $t$ only depends on $t$
and on the state $x$ for the controlled system at the same time. This is a natural
and frequently adopted  formulation. Different formulations are
possible, but usually the corresponding value function is the same
and, if an optimal control exists, it is of feedback type.

\end{remark}

\begin{remark}
All the results that follows admit natural extensions to
slightly more general cases. For instance,  $\lambda$ might
 depend  on time, or the set of admissible control actions
may depend on the present state (so
admissible control laws should satisfy $\alpha(t,x)\in A(x)$,
where $A(x)$ is a given subset  of $A$) provided
appropriate measurability conditions are satisfied.
  We limit ourselves
to the previous setting in order to simplify the notation.
\end{remark}

Let us consider
 the Hamilton-Jacobi-Bellman equation (for short, HJB equation) related
 to the optimal control problem: this is
the following nonlinear integro-differential equation on $[0,T]\times E$:
\begin{eqnarray}
-\frac{\partial v}{\partial t}(t,x)&=&\sup_{a \in A}\left(\mathcal{L}^{a}_E v(t,x) + f(t,x,a)\right) ,
\label{HJB}\\
v(T,x)&=& g(x),\label{HJB_term}
\end{eqnarray}
where the operator $\mathcal{L}^{a}_E$ is defined by
\begin{equation}\label{gen_X}
\mathcal{L}^{a}_E \phi(x) = \int_{E}(\phi(y)-\phi(x))\,\lambda(x,a,dy)
\end{equation}
for all $(t,x,a)\in [0,\,T] \times E \times A$ and every bounded Borel-measurable
function $\phi: E \rightarrow \R$.

\begin{definition}\label{def_sol_HJB}
\begin{em}We say that  a Borel-measurable bounded function
$v:[0,T] \times E\to \R$ is a solution
to the HJB equation if
the right-hand side of \eqref{HJB}  is  Borel-measurable and,
for every $x \in E$, \eqref{HJB_term} holds,
  the map $t \mapsto v(t,x)$ is absolutely continuous in $[0,T]$ and  \eqref{HJB} holds
  almost everywhere on $[0,T]$ (the null
  set of points where it possibly fails may  depend on $x$).
\end{em}
\end{definition}

In the analysis of the HJB equation and the control problem
 we will use the following function spaces, defined
 for any metric space $S$:
\begin{enumerate}
\item
 $C_b(S)=\{\phi:S\to\R\; \text{continuous and bounded}\}$,
  \item
 $LSC_b(S)=\{\phi:S\to\R\; \text{lower semi-continuous and bounded}\}$.
 \item
 $USC_b(S)=\{\phi:S\to\R\; \text{upper semi-continuous and bounded}\}$.
\end{enumerate}
 $C_b(S)$, equipped with the supremum norm $\|\phi\|_\infty$, is a Banach space.
 $LSC_b(S)$ and $USC_b(S)$ are closed subsets of $C_b(S)$, hence  complete metric
 spaces with the induced distance.

In the sequel we need the following classical selection theorem. For a proof
we refer for instance to \cite{BS}, Propositions 7.33 and 7.34,
where a more general statement can also be found.

\begin{proposition}\label{prop_selector}
Let $U$ be a metric space, $V$ a metric separable space.
For $F:U\times V\to \R$ set
$$F^*(u)= \sup_{v\in V}F(u,v),\qquad u\in U.
$$
\begin{enumerate}
\item
If $F\in USC_b(U\times V)$ and $V$ is compact then
$F^*\in USC_b(U)$ and there
 exists a Borel-measurable $\phi  :U\to V$
such that
$$
F(u,\phi(u)) = F^*(u) ,
\qquad u\in U.
$$

\item
If $F\in LSC_b(U\times V)$ then
 $F^*\in LSC_b(U)$ and for
every $\epsilon >0$ there exists a Borel-measurable $\phi_\epsilon :U\to V$
such that
$$
F(u,\phi_\epsilon(u)) \ge F^*(u) -\epsilon,
\qquad u\in U.
$$

\end{enumerate}
\end{proposition}

Next we present a well-posedness result
and a verification theorem
for  the HJB equation in the space
$ LSC_b([0,T] \times E)$,
 Theorems \ref{thm_HJB_ben_posta} and \ref{thm_verifica} below.
The use of lower semi-continuous bounded  functions was already commented
in the introduction and will be useful
for the results in section \ref{FeynmanKac}.
A small variation of our arguments also yields corresponding
results  in the class of upper semi-continuous  functions,
which are more natural when dealing with a maximization problem,
see Theorems \ref{thm_HJB_ben_postaUSC} and
\ref{thm_verificaUSC}
that slightly generalize classical results.
We first formulate the
 assumptions we need.
\begin{equation}
\label{Feller}
 \lambda  \text{ is a Feller transition kernel.}
\end{equation}
We recall that this means that
  for every $\phi\in C_b(E)$
 the function $
(x,a)\to \int_E \phi(y)\,\lambda(x,a,dy)$
  is continuous (hence it
belongs to $ C_b( E\times A)$  by \eqref{lambda_unif_finite}).

Next we will assume either that
\begin{equation}
 \label{fgLSC}
f\in LSC_b([0,T]\times E\times A), \quad g\in LSC_b(E),
\end{equation}
or
\begin{equation}
 \label{fgUSC}
f\in USC_b([0,T]\times E\times A), \quad g\in USC_b(E)
\text{ and }
A\text{ is a compact metric space}.
\end{equation}

\begin{theorem}\label{thm_HJB_ben_posta}
Under the assumptions \eqref{lambda_unif_finite}, \eqref{Feller}, \eqref{fgLSC}
there exists a unique  solution $v\in LSC_b([0,T] \times E)$  to the HJB equation
 (in the sense of
Definition \ref{def_sol_HJB}).
\end{theorem}
\proof
We first make a change of unknown function setting $\tilde v(t,x)=e^{-\Lambda t}
v(t,x)$,
where $\Lambda:=\sup_{x \in E, a \in A}\, \lambda(x,a,E) $
is finite by \eqref{lambda_unif_finite}. It is immediate to check
that $v$ is a solution to \eqref{HJB}-\eqref{HJB_term} if and only if
$\tilde v$ is a solution to
\begin{eqnarray}
-\frac{\partial \tilde v}{\partial t}(t,x)
&=&\sup_{a \in A}\left(\mathcal{L}^{a}_E \tilde v(t,x) +e^{-\Lambda t} f(t,x,a)
+\Lambda \tilde v(t,x)
\right) \nonumber
\\
&=&\sup_{a \in A}\left(\int_{E}\tilde v(t,y)\,\lambda(x,a,dy)+
(\Lambda- \lambda(x,a,E))  \tilde v(t,x) +e^{-\Lambda t} f(t,x,a)
\right), \label{HJBmodif}
\\
\tilde{v}(T,x)&=& e^{-\Lambda T}g(x).\label{HJB_term_modif}
\end{eqnarray}
The notion of solution we adopt for
\eqref{HJBmodif}-\eqref{HJB_term_modif} is completely analogous to Definition \ref{def_sol_HJB}
and need not be repeated.  We set $\Gamma_{\tilde v}(t,x) :=\int_t^T\sup_{a \in A}
\gamma_{\tilde v}(s,x,a)\,ds$ where
\begin{equation}\label{gammavtilde}
    \gamma_{\tilde v}(t,x,a):=
\int_{E}\tilde v(t,y)\,\lambda(x,a,dy)+
(\Lambda- \lambda(x,a,E))  \tilde v(t,x) +e^{-\Lambda t} f(t,x,a)
\end{equation}
and note that solving \eqref{HJBmodif}-\eqref{HJB_term_modif} is equivalent
to finding $\tilde v\in LSC_b([0,T]\times E)$ satisfying
$$
{\tilde v}(t,x)=g(x)+
\Gamma_{\tilde v}(t,x), \qquad t\in [0,T],\,x\in E.
$$
We will prove that  $\tilde v\mapsto g+\Gamma_{\tilde v}$
is a well defined
map of $LSC_b([0,T]\times E)$ into itself and it has a unique fixed point,
which is therefore the required solution.

Fix $\tilde v\in LSC_b([0,T]\times E)$. It follows easily from
\eqref{lambda_unif_finite}
  that $\gamma_{\tilde v}$ is bounded and,
if $\sup_{a\in A}\gamma_{\tilde v}(\cdot,\cdot,a)$ is Borel-measurable,
$\Gamma_{\tilde v}$ is bounded as well. Next we prove that
 $\gamma_{\tilde v}$  and
$\Gamma_{\tilde v}$ are lower semi-continuous.
Note that $(x,a)\mapsto \Lambda- \lambda(x,a,E)$ continuous and
\emph{nonnegative} (this is the reason why we introduced
the equation for $\tilde v$), so
$$ (t,x,a)\mapsto
(\Lambda- \lambda(x,a,E))  \tilde v(t,x) +e^{-\Lambda t} f(t,x,a)
$$
is in $LSC_b([0,T]\times E\times A)$.
Since $\lambda$ is Feller, it is known that the map
 \begin{equation}\label{semicontintegrale}
    (t,x,a)\mapsto\int_{E}\tilde v(t,y)\,\lambda(x,a,dy)
 \end{equation}
 is continuous
when  $\tilde v\in C_b([0,T]\times E)$ (see \cite{BS}, Proposition 7.30).
For general $\tilde v\in LSC_b([0,T]\times E)$, there exists a uniformly
bounded and increasing
sequence $\tilde v_n\in C_b([0,T]\times E)$ such that $\tilde v_n\to \tilde v$
pointwise (see \cite{BS}, Lemma 7.14). From the Fatou Lemma we deduce that the map
\eqref{semicontintegrale}
is in $LSC_b([0,T]\times E\times A)$ and
we conclude that  $\gamma_{\tilde v}\in LSC_b([0,T]\times E\times A)$ as well.
Therefore $\sup_{a\in A}\gamma_{\tilde v}(\cdot,\cdot,a)$,
which equals the right-hand side of \eqref{HJBmodif}, is
lower semi-continuous and hence
Borel-measurable. To prove lower semi-continuity of $\Gamma_{\tilde v}$
suppose $(t_n,x_n)\to (t,x)$; then
\begin{eqnarray*}
  \Gamma_{\tilde v}(t_n,x_n) -\Gamma_{\tilde v}(t,x) &=&
  \int_{t_n}^t \sup_{a\in A}\gamma_{\tilde v}(s,x_n,a)\,ds
+ \int_{t}^T (\sup_{a\in A}\gamma_{\tilde v}(s,x_n,a)-
\sup_{a\in A}\gamma_{\tilde v}(s,x,a))\,ds \\
    &\ge &
    - |t-t_n|\,\|\gamma_{\tilde v}\|_\infty
+ \int_{t}^T (\sup_{a\in A}\gamma_{\tilde v}(s,x_n,a)-
\sup_{a\in A}\gamma_{\tilde v}(s,x,a))\,ds  .
\end{eqnarray*}
By the Fatou Lemma
$$
\liminf_{n\to\infty }\Gamma_{\tilde v}(t_n,x_n) -\Gamma_{\tilde v}(t,x)
\ge
  \int_{t}^T \liminf_{n\to\infty } (\sup_{a\in A}\gamma_{\tilde v}(s,x_n,a)-
\sup_{a\in A}\gamma_{\tilde v}(s,x,a))\,ds  \ge 0,
$$
where in the last inequality we have used the lower semi-continuity
of
$\sup_{a\in A}\gamma_{\tilde v}(\cdot,\cdot,a)$.

Since we assume that $g\in LSC_b(E)$
we have thus checked that $\tilde v\mapsto g+\Gamma_{\tilde v}$
maps $LSC_b([0,T]\times E)$ into itself.
 To prove that it has a unique fixed point we note the easy
 estimate based on \eqref{lambda_unif_finite}, valid
 for every $\tilde v',\tilde v'' \in LSC_b([0,T]\times E)$:
 $$
\begin{array}{l}\dis
  \left|
  \sup_{a \in A} \gamma_{\tilde v'}(t,x,a)-
   \sup_{a \in A} \gamma_{\tilde v''}(t,x,a)
\right|
\le
  \sup_{a \in A} \left|\gamma_{\tilde v'}(t,x,a)-
    \gamma_{\tilde v''}(t,x,a)
\right|
  \\
  \dis
  \qquad\le
  \sup_{a \in A}\left(
  \int_{E}|\tilde v'(t,y)-\tilde v''(t,y)|\,\lambda(x,a,dy)
  +  |\tilde v'(t,x)-\tilde v''(t,x)|\,\lambda(x,a,E)
  \right)
   \\
  \dis
  \qquad\le 2\Lambda \, \|\tilde v'-\tilde v''\|_\infty.
\end{array}
 $$
 By a standard technique one proves that a suitable iteration
 of the map $\tilde v\mapsto g+\Gamma_{\tilde v}$
 is a contraction with respect to the distance induced by the supremum norm,
 and hence that map has a unique fixed point.
\endproof

\begin{theorem}\label{thm_HJB_ben_postaUSC}
Under the assumptions \eqref{lambda_unif_finite}, \eqref{Feller}, \eqref{fgUSC}
there exists a unique  solution $v\in USC_b([0,T] \times E)$  to the HJB equation.
\end{theorem}

\proof
The proof is almost the same as in the previous Theorem, replacing
$LSC_b$ with $USC_b$
with obvious
changes. We introduce
$\tilde v$, $\gamma_{\tilde v}$ and $\Gamma_{\tilde v}$
as before and we prove in particular that
$\gamma_{\tilde v}
\in USC_b([0,T]\times E\times A)$.
The only difference is that we can not immediately conclude that
$\sup_{a\in A}\gamma_{\tilde v}(\cdot,\cdot,a)$  is
upper semi-continuous as well.
However, at this point we can apply point 1 of Proposition
\ref{prop_selector} choosing $U=[0,T]\times E$, $V=A$ and
$F= \gamma_{\tilde v}$ and we deduce
that in fact $\sup_{a\in A}\gamma_{\tilde v}(\cdot,\cdot,a)
\in USC_b([0,T]\times E)$. The rest of the proof is the same.
\endproof

\begin{corollary}\label{thm_HJB_ben_postaC}
Under the assumptions \eqref{lambda_unif_finite}, \eqref{Feller},
if $f\in  C_b([0,T]\times E\times A)$, $g\in  C_b(E)$ and $A$ is a compact
metric space then
the solution $v$   to the HJB equation belongs to
  $C_b([0,T] \times E)$.
\end{corollary}

The Corollary follows immediately from the two previous results.
We proceed to   a verification theorem for the
HJB equation.

\begin{theorem}\label{thm_verifica}
Under the assumptions \eqref{lambda_unif_finite}, \eqref{Feller},
\eqref{fgLSC}
the unique  solution $v\in LSC_b([0,T] \times E)$  to the HJB equation
 coincides with the value function $V$.
\end{theorem}

\proof
Let us fix $(t,x)\in [0,T] \times E $.
As in the proof of
Proposition \ref{prop_construction_X_Markov_pure_jump_1} we have the identity
$$
g(X_T)- v(t,X_t) =   \int_{t}^T  \frac{\partial v}{\partial r}(r,X_r)\,dr
+  \int_{(t,T]}\int_{E} (v(r,y) - v(r,X_{r-}))\,p(dr\,dy),
$$
which follows from the absolute continuity of $t\mapsto v(t,x)$,
taking into account that $X$ is constant among
jump times and using the definition of the random measure $p$.
Given an arbitrary
admissible control $\alpha \in \mathcal{A}_{ad}$ we take the expectation
with respect to the corresponding probability
$\mathbb{P}^{t,x}_{\alpha}$. Recalling that
 the    compensator under
$\mathbb{P}^{t,x}$
is $1_{[t,\infty)}(s)\lambda(X_{s-},\,\alpha(s,X_{s-}),\,dy) \, ds$ we obtain
\begin{eqnarray*}
  \mathbb{E}^{t,x}_{\alpha}[g(X_T)]- v(t,X_t)
   &=& \int_{t}^T  \frac{\partial v}{\partial r}(r,X_r)\,dr
   \\&&
+  \int_{(t,T]}\int_{E} (v(r,y) - v(r,X_{r-}))\,
\lambda(X_{r-},\,\alpha(r,X_{r-}),\,dy) \, dr
 \\
  &=& \int_{t}^T  \left(\frac{\partial v}{\partial r}(r,X_r)
  +\mathcal{L}^{\alpha(r,X_{r})}_E v(r,X_r)
  \right)\,dr .
\end{eqnarray*}
Adding $\mathbb{E}^{t,x}_{\alpha}\int_{t}^{T} f(r,\,X_{r}, \alpha(r,X_r))\, dr$
to both sides and rearranging terms we obtain
\begin{equation}\label{relfondamentale}
v(t,x)=
    J(t,x,\alpha)- \mathbb{E}^{t,x}_{\alpha}
    \int_{t}^T  \left\{\frac{\partial v}{\partial r}(r,X_r)
  +\mathcal{L}^{\alpha(r,X_{r})}_E v(r,X_r)+f(r,\,X_{r}, \alpha(r,X_r))
  \right\}\,dr .
\end{equation}
Recalling the HJB equation
and taking into account that $X$ has piecewise constant trajectories
we conclude that the term in curly brackets
$\left\{\ldots \right\}$ is nonpositive and therefore
 we have $v(t,x)\ge J(t,x,\alpha)$ for every admissible
 control.

 Now we recall that in the proof of
 Theorem \ref{thm_HJB_ben_posta}
 we showed  that  the function $\gamma_{\tilde v}$ defined in
  \eqref{gammavtilde} belongs to
$LSC_b([0,T] \times E\times A)$. Therefore the function
$$
F(t,x,a):=
e^{\Lambda t}
    \gamma_{\tilde v}(t,x,a)=
    \mathcal{L}^{a}_E  v(t,x) + f(t,x,a)
+\Lambda   v(t,x)
$$
is also lower semi-continuous and bounded.
Applying  point 2 of
Proposition \ref{prop_selector} with
$U=[0,T]\times E$ and $V=A$
we see that
 for
every $\epsilon >0$ there exists a Borel-measurable
$\alpha_\epsilon :[0,T]\times E\to A$
such that
$
F(t,x,\alpha_\epsilon(t,x)) \ge \inf_{a\in A}F(t,x,a) -\epsilon$
for all $t\in [0,T]$, $x\in E$.
Taking into account the HJB equation we conclude that
for every $x\in E$ we have
$$
\mathcal{L}^{\alpha_\epsilon(t,x)}_E  v(t,x) + f(t,x,\alpha_\epsilon(t,x))
\ge
-\frac{\partial v}{\partial t}(t,x)
-\epsilon
$$
for almost all $t\in[0,T]$. Noting that
$\alpha_\epsilon$ is an admissible control and choosing $\alpha=\alpha_\epsilon$
in \eqref{relfondamentale} we obtain
$v(t,x)\le J(t,x,\alpha_\epsilon)+\epsilon (T-t)$.
Since we know that
$v(t,x)\ge J(t,x,\alpha)$ for every $\alpha\in\cala_{ad}$
we conclude that $v$ coincides with the value function $V$.
\endproof

\begin{theorem}\label{thm_verificaUSC}
Under assumptions \eqref{lambda_unif_finite}, \eqref{Feller},
\eqref{fgUSC}
the unique  solution $v\in USC_b([0,T] \times E)$  to the HJB equation
 coincides with the value function $V$. Moreover there exists an optimal
 control $\alpha$, which is given by any function satisfying
 \begin{equation}\label{supinHJB}
  \mathcal{L}^{\alpha(t,x)}_E v(t,x) + f(t,x,\alpha(t,x))=
    \sup_{a \in A}\left(\mathcal{L}^{a}_E v(t,x) + f(t,x,a)\right) .
 \end{equation}
\end{theorem}

\proof
We proceed as in the previous proof, but we can now
apply   point 2 of
Proposition \ref{prop_selector} to the function $F$ and
deduce that  there exists a Borel-measurable
$\alpha  :[0,T]\times E\to A$
such that
\eqref{supinHJB} holds.
Any such control $\alpha$ is optimal: in fact we obtain
for every $x\in E$,
$$
\mathcal{L}^{\alpha(t,x)}_E  v(t,x) + f(t,x,\alpha(t,x))
=
-\frac{\partial v}{\partial t}(t,x)
$$
for almost all $t\in[0,T]$ and so
$v(t,x)= J(t,x,\alpha)$.
\endproof

\begin{remark}\label{confrontoPliska}
As already mentioned,
Theorems \ref{thm_HJB_ben_postaUSC} and
\ref{thm_verificaUSC} are similar to
 classical results: compare for instance
\cite{Pl}, Theorems 10, 12, 13, 14.
In that paper the author solves the HJB equations by means of
a general result on nonlinear semigroups of operators, and
for this he requires some more functional-analytic structure,
for instance he embeds the set of decision rules into a properly
chosen topological vector space. He also has more stringent conditions
of the kernel $\lambda$, for instance $\lambda(x,a,B)$ should be strictly
positive and continuous in $(x,a)$ for each fixed $B\in \cale$.
\end{remark}

\section{Control randomization and  dual optimal control problem}\label{Section_control_rand}

In this section we start to implement the
 control randomization
method. In the first step,
 for any initial time $t\ge 0$ and starting point $x\in E$, we construct
  an (uncontrolled) Markovian pair of pure jump
 stochastic processes $(X, I)$ with values in $E\times A$,
by specifying its rate transition measure $\Lambda$ as in \eqref{transcoppia} below.
Next we formulate an auxiliary
optimal control problem where, rougly speaking,
we optimize a cost functional by modifying the intensity of the process
$(X, I)$ over a suitable family. This ``dual'' control problem
will be studied in section \ref{Section_BSDE_rep}  by an approach based on BSDEs.
In section \ref{FeynmanKac} we will prove that the dual value function coincides
with the one introduced in the previous section.

\subsection{A randomized control system}\label{X_I_construction}

Let $E$, $A$ be  Borel spaces with corresponding Borel $\sigma$-algebras
  $\cale$, $\cala$ and let $\lambda$ be
a measure transition kernel   from $(E \times A, \cale\otimes \cala)$
to $(E, \mathcal{E})$ as before. As another basic datum we
suppose we are given a  finite measure $\lambda_0$ on $(A, \mathcal{A})$
with full topological support, i.e., it is strictly positive on any non-empty open
subset of $A$. Note that since $A$ is metric separable such a measure
can always be constructed, for instance supported on a dense discrete subset of $A$.
We still assume  \eqref{lambda_unif_finite}, so we formulate
the following assumption:
\medskip

\noindent \textbf{(H$\lambda$)} \qquad
$\lambda_0$ is a finite measure on $(A, \mathcal{A})$ with full topological support
and  $\lambda$ satisfies
\begin{equation}
\sup_{x \in E, a \in A}\, \lambda(x,a,E) < \infty.\label{lambda_transition_kernel}
\end{equation}

We wish to construct a Markov process as in
section  \ref{Section_construction_X}, but with state space
 $E\times A$. Accordingly, let $\Omega'$ denote the set of sequences
 $\omega'=(t_n,e_n,a_n)_{n \geq 1}$
 contained in $ ((0,\infty)\times E\times A )\cup \{(\infty,\Delta,\Delta')\}$, where
 $\Delta\notin E$ (respectively, $\Delta'\notin A$)
 is   adjoined to  $E$  (respectively, to $A$) as an isolated point,
 satisfying \eqref{processo_punto}
In  the sample space
$\Omega=E\times A\times \Omega'$
we define   $T_n:\Omega\to (0,\infty]$,
$E_n:\Omega\to E\cup\{\Delta\}$, $A_n:\Omega\to A\cup\{\Delta'\}$,
as follows:  writing $\omega=(e,a,\omega')$ in the form
$\omega=(e,a,t_1,e_1,t_2,e_2,\ldots)$ we set  for $t\ge 0$ and for $n\ge 1$
$$\begin{array}{l}\dis
T_n (\omega )= t_n ,
\qquad  T_\infty  (\omega )= \lim_{n\to\infty} t_n,
\qquad T_0(\omega)=0,
\\
  E_n(\omega)=e_n,\qquad A_n(\omega)=a_n,
 \qquad E_0(\omega)=e, \qquad A_0(\omega)=a.
\end{array}
 $$
We also define processes $X:\Omega\times [0,\infty)
\to E\cup\{\Delta\}$, $I:\Omega\times [0,\infty)
\to A\cup\{\Delta'\}$ setting
 $$
X_t = 1_{[0, T_1]}(t)\,E_0+\sum_{n\ge 1}
  1_{ (T_n,T_{n+1}]  }(t)\,E_n,
  \qquad
I_t = 1_{[0, T_1]}(t)\,A_0+\sum_{n\ge 1}
  1_{ (T_n,T_{n+1}]  }(t)\,A_n,
$$
for $t<T_\infty$, $X_t=\Delta$ and $I_t=\Delta'$ for $t\ge T_\infty$.

In $\Omega$ we introduce for all $t\ge0$ the $\sigma$-algebras
 $\calg_t=\sigma(N(s,B)\,:\, s\in (0,t],B\in\cale\otimes \cala)$  generated
 by the counting processes $N(s,B)=
 \sum_{n\ge 1}1_{T_n\le s}1_{(E_n,A_n)\in B}$
and  the $\sigma$-algebra $\calf_t$
  generated by $\calf_0$ and $\calg_t$,
where
$\calf_0:=\cale\otimes\cala\otimes \{\emptyset,\Omega'\}$.
We still denote $\F=(\calf_t)_{t\ge 0}$
and  $\calp$  the corresponding filtration and
predictable $\sigma$-algebra.
By abuse
of notation we also denote by the same symbol
 the trace of $\mathcal{P}$ on
  subsets of the form $[0,\,T] \times \Omega$ or  $[t,T] \times \Omega$,
  for deterministic times   $0\le t\le T < \infty$.

 The  random measure $p$ is now  defined on $(0,\infty)\times E\times A$ as
\begin{eqnarray}\label{misura_aleatoria_X_I}
p(ds\, dy\,db)= \sum_{n\in \N}\one_{\{T_n < \infty\}}\,
\delta_{\{T_n,E_n ,A_n \}}(ds\,dy\,db).
\end{eqnarray}

By means of  $\lambda$ and
$\lambda_0$  satisfying assumption
\textbf{(H$\lambda$)}
we define a (time-independent) rate transition measure on $E\times A$ given by
\begin{equation}\label{transcoppia}
\Lambda(x,a; dy\,db) = \lambda(x,a,dy)\,\delta_{a}(db)+\lambda_0(db)\,\delta_{x}(dy) .
\end{equation}
and the corresponding generator $\mathcal{L}$:
\begin{eqnarray*}
\mathcal{L}\varphi(x,a) &:=& \int_{E \times A}(\varphi(y,b) - \varphi(x,a))\, \Lambda(x,a;dy \, db)
\label{generatoredoppio}\\
&=&  \int_{E}(\varphi(y,a) - \varphi(x,a))\, \lambda(x,a,dy)  \,+ \,
\int_{A }(\varphi(x,b) - \varphi(x,a))\, \lambda_0(db),
\end{eqnarray*}
for all $(x,a)\in E \times A$
and every function $\varphi\in B_b( E \times A)$.

Given any starting time $t\ge0$ and starting point $(x,a)\in E\times A$, an
application of
Proposition \ref{prop_construction_X_Markov_pure_jump_1} provides
  a  probability measure on $(\Omega, \mathcal{F}_\infty)$,
denoted by $\mathbb{P}^{t,x,a}$, such that
$(X,\,I)$ is a Markov process on the time interval $[t,\infty)$ with
respect to $\mathbb{F}$ with transition probabilities associated to $\call$.
Moreover, $\mathbb{P}^{t,x,a}$-a.s.,
$X_s=x$ and $I_s = a$ for all $s\in [0,t]$.
Finally, the restriction   of the measure
$p$ to $(t,\infty)\times E\times A$ admits as
 $\F$-compensator under $\mathbb{P}^{t,x,a}$   the random measure
$$
\tilde{p}(ds\,dy\,db):=   \lambda_0 (db) \, \delta_{\{ X_{s-} \}}(dy) \,ds  +
\lambda(X_{s-},\,I_{s-},\,dy)\, \delta_{\{ I_{s-} \}}(db)\, ds.
$$
We denote
$q:= p- \tilde{p}$   the compensated martingale measure associated to $p$.

\begin{remark}
\label{lambda_0_e_salti}
Note that $\Lambda(x,a; \{x,a\})=\lambda_0(\{a\})+\lambda(x,a,\{x\})$.
So even if we assumed that
$\lambda(x,a,\{x\})=0$,
 in general the rate measure $\Lambda$ would not satisfy the corresponding condition
$\Lambda(x,a; \{x,a\})=0$. We remark that imposing the additional requirement
that $\lambda_0(\{a\})=0$ is too restrictive since, due to the assumption
that $\lambda_0$ has full support, it would
rule out the important case when the space of control actions $A$ is
finite or countable.
\end{remark}

\subsection{The dual optimal control problem}

We introduce  a dual control problem  associated to the process $(X,I)$ and
 formulated in a weak form. For fixed $(t,x,a)$,
it consists in  defining
 a   family of probability
 measures
 $\{\P_\nu^{t,x,a},\,\nu \in \mathcal{V}\}$ in the space $(\Omega, \mathcal{F}_\infty)$,
 all absolutely continuous with respect to $\mathbb{P}^{t,x,a}$,
 whose effect is to change the stochastic intensity
 of the process $(X,I)$ (more precisely, under each $\P_\nu^{t,x,a}$
 the compensator of the associated point process takes
 a desired form), with the aim of maximizing a
 cost   depending on  $f,g$.
 We note that $\{\P_\nu^{t,x,a},\,\nu \in \mathcal{V}\}$ is
a dominated family of probability measures. We proceed with precise
definitions.

We still assume that \textbf{(H$\lambda$)} holds.
Let us define
\begin{eqnarray*}
\mathcal{V} &=& \{ \nu:\Omega\times [0,\infty)\times  A\to (0,\infty),\;
\calp\otimes\cala
\text{-measurable and bounded}  \}.
\end{eqnarray*}
For every $\nu \in \mathcal{V}$, we  consider the predictable random measure
\begin{equation}\label{dual_comp}
\tilde{p}^\nu(ds\,dy\,db):= \nu_s(b) \, \lambda_0 (db) \, \delta_{\{ X_{s-} \}}(dy) \, ds +
\lambda(X_{s-},\,I_{s-},\,dy)\, \delta_{\{ I_{s-} \}}(db)\, ds.
\end{equation}
Now we fix $t\in [0,T]$, $x\in E$, $ a \in A$ and, with the help of a theorem of Girsanov type,
 we will show how to construct a probability measure
on $(\Omega, \mathcal{F}_\infty)$, equivalent to $\mathbb{P}^{t,x,a}$,
under which $\tilde p^\nu$ is the compensator of the measure $p$
on $(0,T]\times E\times A$.
 By the Radon-Nikodym theorem
  one can  find two nonnegative functions $d_1$, $d_2$
defined on $\Omega\times [0,\,\infty) \times E \times A$,
measurable with respect to $\calp\otimes \cale\otimes \cala$ such that
\begin{eqnarray*}
\lambda_0(db) \, \delta_{\{ X_{t-} \}}(dy) \, dt &=& d_1(t,y,b)\, \tilde{p}(dt\,dy\,db)
\\
\lambda(X_{t-},\,I_{t-},\,dy) \, \delta_{\{ I_{t-} \}}(db) \, dt &=&
d_2(t,y,b)\, \tilde{p}(dt\,dy\,db),
\\
d_1(t,y,b) + d_2(t,y,b) &=& 1, \qquad \tilde{p}(dt\,dy\,db)-a.e.
\end{eqnarray*}
 and we have
$
d \tilde{p}^{\nu} = (\nu \, d_1  + d_2)\, d\tilde{p}$.
For any $\nu \in \mathcal{V}$, consider then the Doléans-Dade exponential local martingale
$L^\nu$ defined setting $L_s^\nu=1$ for $s\in [0,t]$ and
\begin{eqnarray*}
L_s^\nu &=& \exp\bigg(\int_t^s\!\int_{E\times A}\log(\nu_r(b) \, d_1(r,y,b) \, + d_2(r,y,b))\,p(dr \, dy \,db) - \int_{t}^{s}\!\int_{A}(\nu_r(b) - 1)\lambda_0(db)\,dr\bigg)
\\&=& e^{\int_{t}^{s}\int_{A}(1 - \nu_r(b) )\lambda_0(db)\,dr} \prod_{n \geqslant 1: T_{n} \leqslant s}(\nu_{T_{n}}(A_n)\,d_1(T_{n},E_n,A_n) + d_2(T_{n},E_n,A_n))\\
\end{eqnarray*}
for $s\in [t,T]$, where $q= p- \tilde{p}$.
When $L^\nu$ is a true martingale, i.e., $\spertxa{L_T^\nu} = 1$, we can
define  a probability measure $\P^{t,x,a}_{\nu}$ equivalent to $\P^{t,x,a}$ on $(\Omega,\,\mathcal{F}_\infty)$ setting $\P^{t,x,a}_{\nu}(d\omega)=L_T^\nu(\omega)
\P^{t,x,a}(d\omega)$.
By the Girsanov theorem for point processes (\cite{J_a}, Theorem 4.5)
the restriction of the random measure $p$ to $(0,T]\times E\times A$
admits $\tilde{p}^\nu=(\nu\,d_1 + d_2)\,\tilde{p}$
as compensator   under $\P^{t,x,a}_{\nu}$.
We denote by $\mathbb{E}_\nu^{t,x,a}$ the expectation operator under $\P_\nu^{t,x,a}$
and by $q^\nu: = p-\tilde{p}^\nu$ the compensated martingale measure of $p$ under $\P^{t,x,a}_{\nu}$.
The validity of the condition $\spertxa{L_T^\nu} = 1$  under our assumptions,
as well as other useful properties, are proved in the following proposition.

\begin{lemma}\label{lemma_P_nu_martingale}
Let assumption {\bf (H$\lambda$)} hold. Then, for every
$t\in [0,T]$, $x\in E$ and
$\nu \in \mathcal{V}$, under the probability
$\P^{t,x,a}$
the process $L^\nu$ is a   martingale on $[0,T]$
 and $L_T^\nu$ is square integrable.

In addition, for every  $\mathcal{P} \otimes \cale \otimes \cala$-measurable  function $H: \Omega \times [t,T] \times E \times A \rightarrow \R$ such that $\spertxa{\int_t^T \int_{E\times A} |H_s(y,b)|^2\, \tilde{p}(ds\,dy\,db)} < \infty$,  the process $\int_{t}^{\cdot}\int_{E \times A}H_s(y,b)\, q^\nu(ds\,dy\, db)$ is a $\P^{t,x,a}_{\nu}$-martingale on $[t,T]$.
\end{lemma}
\proof
The first part of the proof is inspired by Lemma 4.1 in \cite{KhPh}. In particular, since $\nu$ is bounded and $\lambda_0(A)<\infty$, we see that
\[
S_T^\nu = \exp\bigg(\int_t^T\int_A |\nu_s(b) - 1|^2 \lambda_0(db)\,ds \bigg)
\]
is bounded. Therefore, from Theorem 8, see also Theorem 9, in \cite{Pro-Shi},  follows the martingale property of $L^\nu$ together with its uniform integrability. Concerning the square integrability of $L_T^\nu$, set $\ell(x,\lambda):=2\ln(x\lambda + 1-\lambda)-\ln(x^2\lambda+1-\lambda)$, for any $x\geq0$ and $\lambda\in[0,1]$. From the definition of $L^\nu$ we have (recalling that $d_2(s,y,b)=1-d_1(s,y,b)$)
\[
|L_T^\nu|^2 = L_T^{\nu^2}S_T^\nu \exp\bigg(\int_t^T\int_{E\times A} \ell(\nu_s(b),d_1(s,y,b))\,p(ds\,dy\,db)\bigg) \leq L_T^{\nu^2}S_T^\nu,
\]
where the last inequality follows from the fact that $\ell$ is nonpositive. This entails that $L_T^\nu$ is square integrable.

Let us finally  fix  a predictable  function $H$ such that
$\spertxa{\int_t^T \int_{E\times A} |H_s(y,b)|^2\, \tilde{p}(ds\,dy\,db)} < \infty$.
The process $\int_{t}^{\cdot}\int_{E\times A} H_s(y,b)\, q^{\nu}(ds\,dy\,db)$ is
a $\P^{t,x,a}_{\nu}$-local martingale, and the uniform integrability follows from  the Burkholder-Davis-Gundy and
 Cauchy Schwarz inequalities, together with the square integrability of  $L_T^\nu$.
\endproof

To complete the formulation of the
dual optimal control problem we specify the conditions that we will
assume for the cost functions $f$, $g$:

\medskip

\noindent \textbf{(Hfg)}\qquad $f\in B_b([0,\,T] \times E \times A)$
and  $g \in B_b(E)$.

\medskip

For every
$t\in [0,T]$, $x\in E$, $a \in A$ and
$\nu \in \mathcal{V}$
we finally introduce the dual gain functional
$$
J(t,x,a,\nu)= \mathbb{E}^{t,x,a}_{\nu}\left[g(X_{T}) + \int_{t}^{T} f(s,\,X_{s},\,I_{s}) ds\right],
$$
and the dual value function
\begin{equation}\label{funzione_valore_duale}
    v^{\ast}(t,x,a) = \sup_{ \nu \in \mathcal{V}} J(t,x,a,\nu).
\end{equation}

\begin{remark}
An interpretation of the dual optimal control problem can be given as follows.
It can be proved that, under $\P_\nu^{t,x,a}$,
 the jump times of $I$ and $X$,  denoted by $\{S_n\}$ and $\{R_n\}$ respectively,
   are disjoint. The compensators of the corresponding
random measures
$ \mu^I(ds\,db)= \sum_{n}\, \delta_{(S_n, I_{S_n})}\,(ds\,db)$ on
$(0,\infty)\times A$ and
$ \mu^X(ds\,dy)= \sum_{n}\, \delta_{(R_n, X_{R_n})}\,(ds\,dy)$ on
$(0,\infty)\times E$ are
$$\tilde{\mu}^I(ds\,db)= \nu_s(b)\,\lambda_0(db)\,ds,\qquad
\tilde{\mu}^X(ds\,dy)= \lambda(X_{s-},I_{s-},dy)\,ds.
$$
Thus, the effect of choosing  $\nu$ is to change the intensity
of the $I$-component. We leave the proofs of these facts to
the reader since they will not be used in the sequel.
\end{remark}

\section{Constrained BSDE and representation of the dual value function}
\label{Section_BSDE_rep}

In this section we introduce a BSDE, with a sign constrain on its
martingale part, and prove existence and uniqueness
of a minimal solution, in an appropriate sense. The BSDE
is then used to give a representation formula for the dual value function
introduced above.

Throughout this section we assume that the
assumptions \textbf{\textup{(H$\lambda$)}} and \textbf{\textup{(Hfg)}}
are satisfied
and we use the randomized control setting introduced above:
$\Omega, \F, X, \P^{t,x,a}$ as well as the random measures
$ p, \tilde p, q $
are the same as in in subsection
\ref{X_I_construction}.
For any $(t,x,a) \in [0,\,T] \times E \times A$, we introduce the following notation.
\begin{itemize}
    \item $\mathbb{\textbf{L}}^{\textbf{2}}(\lambda_0)$,  the set of $\mathcal{A}$-measurable maps $\psi:  A \rightarrow \R$ such that
    \begin{eqnarray*}
    |\psi|^2_{\mathbb{\textbf{L}}^{\textbf{2}}(\lambda_0)} &:=& \int_{A} |\psi(b)|^2 \, \lambda_0(db) < \infty.
    \end{eqnarray*}
\item $\mathbb{\textbf{L}}^{\textbf{2}}_{\textbf{t,x,a}}(\mathcal{F}_\tau)$,
 the set of $\mathcal{F}_\tau$-measurable random variable $X$ such that $\spertxa{|X|^2} < \infty$; here $\tau$ is an $\F$-stopping time with values in
    $[t,T]$.
\item $\mathbb{\textbf{S}}^{\textbf{2}}_{\textbf{t,x,a}}$ the set of real valued càdlàg adapted  processes $Y = (Y_s)_{t \leqslant s \leqslant T}$ such that
    \begin{displaymath}
    ||Y||^2_{\mathbb{\textbf{S}}^{\textbf{2}}_{\textbf{t,x,a}}} := \spertxa{\sup_{t \leqslant s \leqslant T}|Y_s|^2} < \infty.
    \end{displaymath}
\item $\mathbb{\textbf{L}}^{\textbf{2}}_{\textbf{t,x,a}}\textup{(q)}$,  the set of $\mathcal{P} \otimes \mathcal{E} \otimes \mathcal{A}$-measurable maps $Z: \Omega \times [t,T] \times E \times A \rightarrow \R$ such that
   $$
    \begin{array}{l}\dis
    ||Z||^2_{\mathbb{\textbf{L}}^{\textbf{2}}_{\textbf{t,x,a}}(q)} := \spertxa{ \int_{t}^{T}\int_{E \times A} |Z_s(y,b)|^2 \, \tilde{p}(ds\,dy\,db)   }\\
    \dis\quad
    = \spertxa{   \int_{t}^{T}\int_{E} |Z_s(I_s,y)|^2 \, \lambda(X_s,I_s,dy)\,ds + \int_{t}^{T}\int_{A} |Z_s(X_s,b)|^2 \, \lambda_0(db)\,ds   } < \infty.
    \end{array}
    $$
\item $\mathbb{\textbf{K}}^{\textbf{2}}_{\textbf{t,x,a}}$ the set of nondecreasing predictable processes $K = (K_s)_{t \leqslant s \leqslant T} \in \mathbb{\textbf{S}}^{\textbf{2}}_{\textbf{t,x,a}}$ with $K_t = 0$, with the induced
    norm
    \begin{displaymath}
    ||K||^2_{\mathbb{\textbf{K}}^{\textbf{2}}_{\textbf{t,x,a}}} = \spertxa{|K_T|^2}.
    \end{displaymath}
\end{itemize}
We are interested in studying the following family of BSDEs  parametrized by $(t,x,a)$: $\mathbb{P}^{t,x,a}$-a.s.,
\begin{eqnarray}\label{BSDE}
Y_{s}^{t,x,a} &=& g(X_T) + \int_{s}^{T}f(r,X_r,I_r)\,dr + K_T^{t,x,a} - K_s^{t,x,a}\\
&& - \int_{s}^{T}\int_{E \times A}Z_r^{t,x,a}(y,\,b)\, q(dr\,dy\, db)  - \int_{s}^{T}\int_{A}Z_r^{t,x,a}(X_r,\,b)\, \lambda_0(db)\,dr, \,\, s\in[t,T],\nonumber
\end{eqnarray}
with the sign constraint
\begin{equation}\label{BSDE_constraint}
Z_s^{t,x,a}(X_{s-},b)\leqslant 0, \qquad  ds \otimes d\mathbb{P}^{t,x,a} \otimes \lambda_0(db)-\text{a.e.}\,\,\text{on}\,\, [t,T]\times \Omega \times A.
\end{equation}
This constraint can be seen as a sign condition imposed on the jumps
of the corresponding
stochastic integral.
\begin{definition}
A  solution to the equation \eqref{BSDE}-\eqref{BSDE_constraint} is a triple $(Y,Z,K)\in \mathbb{\textup{\textbf{S}}}^{\textup{\textbf{2}}}_{\textup{\textbf{t,x,a}}} \times \mathbb{\textup{\textbf{L}}}^{\textup{\textbf{2}}}_{\textup{\textbf{t,x,a}}}\textup{(q)}\times \mathbb{\textup{\textbf{K}}}^{\textup{\textbf{2}}}_{\textup{\textbf{t,x,a}}}$
that satisfies  \eqref{BSDE}-\eqref{BSDE_constraint}.

A  solution $(Y,Z,K)$ is called minimal if for any other
 solution  $(\tilde{Y}, \tilde{Z},\tilde{K})$ we  have, $\P^{t,x,a}$-a.s.,
\begin{displaymath}
Y_s \leqslant \tilde{Y}_s, \qquad  s\in[t,T].
\end{displaymath}
\end{definition}
%We assume that the terminal condition $g$ and the generator $f$ satisfy the following conditions.
%
%\medskip
%
%\noindent \textbf{(HBC)}
%\begin{itemize}
%\item[(i)]$g: E \rightarrow \R$ is a $\mathcal{F}^t_T \otimes \mathcal{E}$-measurable  function;
%%such that
%%\begin{displaymath}
%%\spertxa{|g(X_T^{t,x,\alpha})|^2} < \infty;
%%\end{displaymath}
%\item[(ii)] $f: [0,\,T] \times E \times A \rightarrow \R$ is a $\mathcal{P}^t_T \otimes \mathcal{E}\otimes \mathcal{A}$-measurable function;
% %such that
%    % for every $t \in [0,\,T]$, $x \in E$, $a \in A$ and $\alpha(\cdot)\in \mathcal{A}_{\mathbb{F}_0}$,
%%\begin{displaymath}
%%\sup_{a \in A}f(s, x, a) < \infty;
%%\end{displaymath}
%%\begin{displaymath}
%%\spertxa{\int_{t}^{T}\sup_{a \in A}|f(s, X_s, a)|^2 \,dt} < \infty;
%%\end{displaymath}
%%there exists $\beta >1$ such that
%%\begin{displaymath}
%%\spertxa{\int_{t}^{T}|f(s, X_s^{t,x,\alpha}, \alpha_s)| \,ds}^{\beta}< \infty;
%%\end{displaymath}
%\item[(iii)] for every $t \in [0,\,T]$, $x \in E$ and $a \in A$, there exists a nonnegative constant $M$ such that
%\begin{equation}\label{bound_f_g}
%|f(t,x,a)| + |g(x)| \leqslant M.
%\end{equation}
%\end{itemize}

\begin{proposition}\label{Prop_uniq_min_sol}
Under assumptions \textbf{\textup{(H$\lambda$)}} and \textbf{\textup{(Hfg)}}, for any $(t,x,a)\in [0,\,T] \times E \times A$, if there exists a minimal solution on $(\Omega, \mathcal{F}, \mathbb{F}, \P^{t,x,a})$ to the BSDE \eqref{BSDE}-\eqref{BSDE_constraint}, then it is unique.
\end{proposition}
\proof
Let $(Y,Z,K)$ and $(Y',Z',K')$ be two minimal solutions of \eqref{BSDE}-\eqref{BSDE_constraint}.
The component $Y$ is unique  by definition, and the difference between the two backward equations gives: $\P^{t,x,a}$-a.s.
$$
\begin{array}{l}\dis
\int_{t}^{s}\int_{E \times A}(Z_r(y,\,b)-Z'_r(y,\,b))\, p(dr\,dy\, db)
\\\qquad= K_s - K'_s
 \dis+ \int_{t}^{s}\int_{E}(Z_r(y,\,I_{r-})-Z'_r(y,\,I_{r-}))\, \lambda(X_{r-},I_{r-}dy)\,dr,\quad \forall \,t \leqslant s \leqslant T.
\end{array}
$$
The right hand is a predictable process, in particular it has no totally inaccessible jumps (see, e.g., Proposition 2.24, Chapter I, in \cite{J-S}), while the left side is a pure jump process with totally inaccessible jumps, unless $Z = Z'$. This implies  the uniqueness of the component $Z$, and as a consequence the component $K$ is unique as well.
\endproof

We now state the main result of the section.
\begin{theorem}\label{Thm_ex_uniq_minimal_BSDE}
Under the assumptions \textbf{\textup{(H$\lambda$)}} and \textbf{\textup{(Hfg)}}, for all $(t,x,a) \in [0,\,T] \times E \times A$ there exists a unique minimal solution $Y^{t,x,a}$ to \eqref{BSDE}-\eqref{BSDE_constraint}. Moreover, for all $s \in [t,T]$, $Y_s^{t,x,a}$  has the explicit representation:
\begin{equation}\label{rep_Y}
Y_s^{t,x,a} = \esssup_{\nu \in \mathcal{V}}
\spernu{g(X_{T}) + \int_{s}^{T} f(r,X_{r},I_{r})\, dr \,\bigg\vert\, \calf_s}, \quad s \in [t,T].
\end{equation}
In particular, setting $s=t$, we have
the following representation formula for the value function of the dual control problem:
\begin{equation}\label{vstar_rep}
v^{\ast}(t,x,a)=Y_t^{t,x,a}, \qquad (t,x,a)\in[0,T]\times E\times A.
\end{equation}
\end{theorem}

The rest of this section is devoted to prove Theorem \ref{Thm_ex_uniq_minimal_BSDE}.
To this end we will  use a penalization approach presented in the following subsections.
Here we only note that for the solvability of the BSDE the use of the
filtration $\F$ introduced above is essential, since it involves application of
martingale representation theorems for multivariate
point processes (see e.g. Theorem 5.4 in \cite{J_a}).

\subsection{Penalized BSDE and associated dual control problem}
Let us consider the family of penalized BSDEs associated to \eqref{BSDE}-\eqref{BSDE_constraint}, parametrized by the integer $n \geqslant 1$: $\P^{t,x,a}$-a.s.,
\begin{eqnarray}\label{BSDE_penalized_K}
Y_{s}^{n,t,x,a}&=& g(X_{T}) + \int_{s}^{T} f(r,X_{r},I_{r})\, dr + K_T^{n,t,x,a} - K_s^{n,t,x,a}\\
&& - \int_{s}^{T} \int_{E \times A} Z_{r}^{n,t,x,a}(y,b) \, q(dr\,dy\,db) - \int_{s}^{T} \int_{A} Z_{r}^{n,t,x,a}(X_{r},b) \,\lambda_{0}(db)\, dr,
 \,\,\, s \in [t,T],\nonumber
\end{eqnarray}
where $K^n$ is the nondecreasing process in $\mathbb{\textbf{K}}^{\textbf{2}}_{\textbf{t,x,a}}$ defined by
\begin{displaymath}
K_s^n = n \,\int_{t}^{s} \int_{A}[Z_{r}^n(X_{r},b)]^{+} \,\lambda_{0}(db)\, dr.
\end{displaymath}
Here we denote by $[u]^+$ the positive part of $u$.
The penalized BSDE \eqref{BSDE_penalized_K} can be rewritten in the equivalent form: $\P^{t,x,a}$-a.s.,
\begin{eqnarray*}%\label{BSDE_penalized}
Y_{s}^{n,t,x,a} &=& g(X_{T}) + \int_{s}^{T} f^n(r,\,X_{r},\,I_{r},\,Z_r^{n,t,x,a}(X_r,\cdot))\, ds\\
&&\qquad - \int_{s}^{T} \int_{E \times A} Z_{r}^{n,t,x,a}(y,b) \, q(dr\,dy\,db), \quad s \in [t,T].\nonumber
\end{eqnarray*}
where the generator $f^n$ is defined by
\begin{equation}\label{f^n}
f^n(t,x,a,\psi) := f(t,x,a) + \int_{A}\left\{ n\,[\psi(b)]^{+} - \psi(b)\right\} \,\lambda_{0}(db),
\end{equation}
for all $(t,x,a)$ in $[0,\,T] \times E \times A$, and $\psi\in \mathbb{\textbf{L}}^{\textbf{2}}(\lambda_0)$.
We note that
under \textbf{(H$\lambda$)} and \textbf{(Hfg)}   $f^n$ is Lipschitz continuous
  in $\psi$  with respect to the  norm of $\mathbb{\textbf{L}}^{\textbf{2}}(\lambda_0)$,
  uniformly in $(t,x,a)$, i.e., for every $n \in \N$
 there exists a constant $L_n$ depending only on $n$ such that for every
  $(t,x,a)\in [0,\,T] \times E \times A$
   and
  $\psi,\,\psi' \in \mathbb{\textbf{L}}^{\textbf{2}}(\lambda_0)$,
\begin{displaymath}
|f^n(t,x,a,\psi')-f^n(t,x,a,\psi)| \leqslant L_n\,|\psi-\psi'|_{\mathbb{\textbf{L}}^{\textbf{2}}(\lambda_0)}.
\end{displaymath}
 The use of the natural filtration $\F$ allows to use
well known
integral representation results for $\F$-martingales  (see, e.g., Theorem 5.4 in \cite{J_a})
and we have the following proposition, whose proof is standard
and is therefore omitted (similar proofs can be found
in \cite{Xia} Theorem 3.2, \cite{Bech} Proposition 3.2,
\cite{CoFu-m}
Theorem 3.4).

\begin{proposition}
Let assumptions \textbf{\textup{(H$\lambda$)}} and \textbf{\textup{(Hfg)}} hold. For every initial condition $(t,x,a) \in [0,\,T] \times E \times A$, and for every $n \in \N$, there exists a unique solution $(Y_s^{n,t,x,a},Z_s^{n,t,x,a})_{s \in [t,T]} \in \mathbb{\textbf{\textup{S}}}^{\textbf{\textup{2}}}_{\textbf{\textup{t,x,a}}} \times \mathbb{\textbf{\textup{L}}}^{\textbf{\textup{2}}}_{\textbf{\textup{t,x,a}}}\textup{(q)}$ satisfying the penalized BSDE \eqref{BSDE_penalized_K}.
\end{proposition}

Next we show that
 the solution to the penalized  BSDE \eqref{BSDE_penalized_K}
 provides an explicit representation of the value function
of a corresponding dual control problem depending on $n$.
This is the content of Lemma \ref{Lemma_rep_Y_n}
which will
allow to deduce some estimates uniform with respect to  $n$.

For every $n \geqslant 1$, let $\mathcal{V}^n$ denote the subset of elements $\nu \in \mathcal{V}$ that take values in $(0,\,n]$.

\begin{lemma}\label{Lemma_rep_Y_n}
Let assumptions \textbf{\textup{(H$\lambda$)}} and \textbf{\textup{(Hfg)}} hold.
For all $n \ge 1$ and $s \in [t,T]$,
\begin{equation}\label{rep_Y_n}
Y_s^{n,t,x,a} = \esssup_{\nu \in \mathcal{V}^n}\spernu{g(X_{T})
+ \int_{s}^{T} f(r,X_{r},I_{r})\, dr\,\bigg\vert\,\calf_s}, \qquad \P^{t,x,a}-a.s.
\end{equation}
\end{lemma}
\proof
We fix $n \ge 1$ and for any $\nu \in \mathcal{V}^n$ we introduce the compensated martingale measure $q^\nu(ds \, dy \,db)= q(ds \, dy \,db)- (\nu_s(b)-1)\,d_1(s,y,b)\,\tilde{p}(ds\,dy\,db)$ under $\P_\nu^{t,x,a}$.  We see that the solution $(Y^n, Z^n)$ to the BSDE \eqref{BSDE_penalized_K}  satisfies: $\P^{t,x,a}$-a.s.,
\begin{eqnarray}\label{BSDE_proof}
Y_{s}^n &=&  g(X_{T}) + \int_{s}^{T} f(r,X_{r},I_{r})\, dr
+ \int_{s}^{T} \int_{A} \{n [Z_{r}^n(X_r,b)]^+ - \nu_r(b) \,Z_{r}^n(X_r,b)\}\, \lambda_0(db)\, dr \nonumber\\
&& - \int_{s}^{T} \int_{E \times A} Z_{r}^n(y,b)\, q^\nu(dr\,dy\,db), \quad s \in [t,T].
\end{eqnarray}
By taking conditional expectation in \eqref{BSDE_proof} under $\P_\nu^{t,x,a}$
and applying  Lemma \ref{lemma_P_nu_martingale} we get,
for any $s \in [t,T]$,
\begin{eqnarray}\label{BSDE_spernu}
Y_{s}^{n,t,x,a} &=& \spernu{g(X_{T}) + \int_{s}^{T} f(r,X_{r},I_{r})\, dr
\,\bigg\vert\,\calf_s
}\\
&&+\, \spernu{\int_{s}^{T} \int_{A} \{n [Z_{r}^{n,t,x,a}(X_r,b)]^+ - \nu_r(b)\,Z_{r}^{n,t,x,a}(X_r,b)\}\, \lambda_0(db)\, dr
\,\bigg\vert\,\calf_s
},\nonumber
\end{eqnarray}
$\P_\nu^{t,x,a}$-a.s.
From the elementary numerical inequality: $n[u]^+ - \nu u \geqslant 0$ for all $u \in \R$, $\nu \in (0,n]$, we deduce by \eqref{BSDE_spernu} that
\begin{equation}\label{ineq_esssup}
Y_s^{n,t,x,a} \geqslant \esssup_{\nu \in \mathcal{V}^n} \spernu{g(X_{T}) + \int_{s}^{T} f(r,X_{r},I_{r})\, dr\,\bigg\vert\,\calf_s}.
\end{equation}
On the other hand, for $\epsilon \in (0,\,1)$, let us consider the process $\nu^\epsilon \in \mathcal{V}^n$ defined by
\begin{displaymath}
\nu^\epsilon_s(b) = n\,\one_{\{Z_s^{n,t,x,a}(X_{s-},b) \geqslant 0 \}} + \epsilon\,\one_{\{-1 < Z_s^{n,t,x,a}(X_{s-},b) < 0\}} - \epsilon\,Z_s^{n,t,x,a}(X_{s-},b)^{-1}\,\one_{\{Z_s^{n,t,x,a}(X_{s-},b) \leqslant -1\}}.
\end{displaymath}
By construction, we have
\begin{displaymath}
n [Z_{s}^{n,t,x,a}(X_{s-},b)]^+ - \nu^{\epsilon}_s(b)\,Z_{s}^{n,t,x,a}(X_{s-},b) \leqslant \epsilon, \qquad s \in [t,T],\,b \in A,
\end{displaymath}
and thus for the choice of $\nu = \nu^\epsilon$ in \eqref{BSDE_spernu}:
\begin{eqnarray*}
Y_{s}^{n,t,x,a} &\leqslant&  \spertxanuepsilon{g(X_{T}) + \int_{s}^{T} f(r,X_{r},I_{r})\, dr
\,\bigg\vert\,\calf_s
}
+ \epsilon T |\lambda_0(A)|\\
&\leqslant&  \esssup_{\nu \in \mathcal{V}^n}
\spernu{g(X_{T}) + \int_{s}^{T} f(r,X_{r},I_{r})\, dr\,\bigg\vert\,\calf_s}
+ \epsilon T |\lambda_0(A)|.
\end{eqnarray*}
Together with \eqref{ineq_esssup}, this is enough to prove the required representation of $Y^n$.
Note that we could not take $\nu_s(b) = n \one_{\{Z_s^n(X_{s-},b)\geqslant 0\}}$, since this process does not belong to $\mathcal{V}^n$ because of the requirement of strict positivity.
\endproof

\subsection{Limit behavior of the penalized BSDEs and conclusion of the proof of Theorem \ref{Thm_ex_uniq_minimal_BSDE}}

As a consequence of the representation \eqref{rep_Y_n}
we immediately obtain the following  estimates:
\begin{lemma}\label{Lemma_monotonicity_boundedness_Y_n}
Let  assumptions \textbf{\textup{(H$\lambda$)}} and \textbf{\textup{(Hfg)}} hold.
 There exists a constant $C$, depending only on $T,f,g$, such that
for any $(t,x,a)\in [0,\,T]\times E \times A$
and $n \ge 1$, $\P^{t,x,a}$-a.s.,
\begin{displaymath}
Y_s^{n,t,x,a} \leqslant Y_s^{n+1,t,x,a},
\qquad  |Y_s^{n,t,x,a}|\leqslant C, \qquad   s \in [t,T].
\end{displaymath}
\end{lemma}
\proof
For fixed $s\in [t,T]$,
 the
almost sure
 monotonicity of  $ Y^{n,t,x,a}$   follows from the representation formula \eqref{rep_Y_n}, since by definition $\mathcal{V}^n \subset \mathcal{V}^{n+1}$;
  moreover, the same formula shows that we can take
  $C=||g||_{\infty} + T\,||f||_{\infty}$. Finally, these inequalities hold
  for every $s\in [t,T]$ outside a null set, since the processes $Y^{n,t,x,a}$
   are càdlàg.
\endproof

Moreover, the following a priori uniform estimate on the sequence $(Y^{n,t,x,a},Z^{n,t,x,a},K^{n,t,x,a})$ holds:
\begin{lemma}\label{lemma_BSDE_estimations}
Let assumptions \textbf{\textup{(H$\lambda$)}} and \textbf{\textup{(Hfg)}} hold. For all $(t,x,a)\in [0,\,T]\times E \times A$ and $n \in \N$, there exists a positive constant $C'$ depending only on $T,f,g$
such that
\begin{equation}\label{unif_pen_estimate}
||Y^{n,t,x,a}||^2_{\mathbb{\textbf{\textup{S}}}^{\textbf{\textup{2}}}_{\textbf{\textup{t,x,a}}}}
+||Z^{n,t,x,a}||^2_{\mathbb{\textbf{\textup{L}}}^{\textbf{\textup{2}}}_{\textbf{\textup{t,x,a}}}\textup{(q)}} + ||K^{n,t,x,a}||^2_{\mathbb{\textbf{\textup{K}}}^{\textbf{\textup{2}}}_{\textbf{\textup{t,x,a}}}}
  \leqslant C' .
\end{equation}
\end{lemma}
\proof
In the following we omit for simplicity of notation the dependence on $(t,x,a)$ for the triple $(Y^{n,t,x,a},Z^{n,t,x,a},K^{n,t,x,a})$. The estimate on
$Y^{n}$ follows immediately from the previous lemma:
\begin{equation}\label{Yn_unif_estim}
||Y^{n}||^2_{\mathbb{\textbf{S}}^{\textbf{2}}_{\textbf{t,x,a}}} =\spertxa{\sup_{s \in [t,T]}\, |Y_s^{n}|^2} \leqslant C^2.
\end{equation}
Next we notice that, since $K^n$ is continuous, the jumps of $Y^n$ are given by the formula
$$\Delta Y^n_s = \int_{E \times A}Z_s^{n}(y,b)\,p(\{s\},dy \,db).
$$
The It$\hat{\mbox{o}}$ formula applied to $|Y_t^n|^2$ gives:
\begin{eqnarray}\label{ito_jump}
d|Y_{r}^n|^2 &=& 2\,Y_{r-}^n\, d Y_r^n + |\Delta Y_r^n|^2\nonumber\\
&=& -2\,Y_{r-}^n\, f(X_{r-}, I_{r-})\,dr - 2\,Y_{r-}^n\,d K^n_r\nonumber\\
&&+2\,Y_{r-}^n\,\int_{E \times A} Z_r^n(y,b)\, q(dr\,dy\,db)  +2\,Y_{r-}^n\,\int_{A} Z_r^n(X_{r-},b)\, \lambda_0(db)\,dr\nonumber\\
&&+ \int_{E \times A} |Z_r^n(y,b)|^2\, p(\{r\}\,dy\,db).
\end{eqnarray}
Integrating \eqref{ito_jump} on $[s,\,T]$, for every  $s \in [t,T]$, and recalling the elementary inequality $2ab \leqslant \frac{1}{\delta}a^2 + \delta b^2$ for any  constant $\delta >0$, and that
\begin{equation}\label{compensparz}
    \spertxa{\int_{s}^{T}\int_{A}|Z^n_r(X_{r-},\,b)|^2 \, \lambda_0(db)\,dr} \leqslant \spertxa{\int_{s}^{T}\int_{E \times A}|Z^n_r(y,\,b)|^2 \, \tilde{p}(dr\,dy\, db)},
\end{equation}
we have:
\begin{eqnarray}\label{partial_estim1}
&& \spertxa{|Y_{s}|^2} + \spertxa{\int_{s}^{T}\int_{E \times A}|Z^n_r(y,\,b)|^2 \, \tilde{p}(dr\,dy\, db)}\nonumber\\
&& \qquad \qquad \leqslant \spertxa{|g(X_T)|^2}  \nonumber\\
&&\qquad \qquad \, + \frac{1}{\beta} \spertxa{\int_{s}^{T}|f(r,X_r,I_r)|^2\,dr} + \beta \spertxa{\int_s^T|Y_r^n|^2\,dr}\nonumber\\
&&\qquad \qquad \,+ \frac{T\,\lambda_0(A)}{\gamma} \spertxa{\int_{s}^{T}\int_{E \times A} |Z_r^n(y,b)|^2\,\tilde{p}(dr\,dy\, db)} + \gamma \spertxa{\int_s^T|Y_r^n|^2\,dr}\nonumber\\
&&\qquad \qquad\, + \frac{1}{\alpha}\spertxa{\sup_{s \in [t,T]} |Y_s^n|^2} + \alpha \spertxa{ |K_T^n- K_s^n|^2}, \quad s \in [t,T],
\end{eqnarray}
for some $\alpha, \beta, \gamma >0$,
Now, from the equation \eqref{BSDE_penalized_K}
 we obtain:
\begin{eqnarray*}
K_T^n -K_s^n
&=& Y_s^n - g(X_T) -\int_s^T f(r,X_r,I_r)dr\\
&& + \int_{s}^{T}\int_A \,Z_r^n(X_r,b)\,\lambda_0(db)\,dr\\
&& + \int_{s}^{T}\int_{E \times A} \,Z_r^n(y,b)\,q(dr\,dy\,db),
%&=& Y_s^n - g(X_T) -\int_s^T f(r,X_r,I_r)dr\\
%&& + \int_{s}^{T}\int_{E \times A} \,Z_r^n(y,b)\,p(dr\,dy\,db)\\
%&& - \int_{s}^{T}\int_E \,Z_r^n(y,I_{r})\,\lambda(X_{r},\,I_{r},\,dy)\,dr,
\quad s \in [t,T].
\end{eqnarray*}
Next we note the equality
\begin{eqnarray*}
	\spertxa{\bigg |\int_{s}^{T}\int_{E \times A}Z^n_r(y,\,b) \, q(dr\,dy\, db)\bigg|^2} &=& \spertxa{\int_{s}^{T}\int_{E \times A}|Z^n_r(y,\,b)|^2 \, p(dr\,dy\, db)}\\
	&=& \spertxa{\int_{s}^{T}\int_{E \times A}|Z^n_r(y,\,b)|^2 \, \tilde{p}(dr\,dy\, db)}
\end{eqnarray*}
that can be proved applying the Ito formula as before
to the square of the martingale
$u\mapsto \int_{s}^{u}\int_{E \times A}Z^n_r(y,\,b) \, q(dr\,dy\, db)$, $u\in [s,T]$
(or by considering its quadratic variation). Recalling again \eqref{compensparz} we see that
there exists some  positive constant $B$ such that
\begin{eqnarray}\label{partial_estim_K^n}
\spertxa{|K_T^n -K_s^n|^2} &\leqslant& B\bigg( \spertxa{|Y_s^n|^2} + \spertxa{|g(X_T)|^2} +\spertxa{\int_s^T |f(r,X_r,I_r)|^2 \,dr}\nonumber\\
&& + \spertxa{\int_{s}^{T}\int_{E \times A} \,|Z_r^n(y,b)|^2\,\tilde{p}(dr\,dy\,db)}\bigg),
\quad s \in [t,T].
\end{eqnarray}
Plugging \eqref{partial_estim_K^n} into \eqref{partial_estim1}, and recalling the uniform estimation \eqref{Yn_unif_estim} on $Y^n$, we get
\begin{eqnarray}\label{partial_estim}
&& (1 - \alpha B)\,\spertxa{|Y_{s}|^2} + \left( 1 - \left[ \alpha B + \frac{T\,\lambda_0(A)}{\gamma}\right] \right)\,\spertxa{\int_{s}^{T}\int_{E \times A}|Z^n_r(y,\,b)|^2 \, \tilde{p}(dr\,dy\, db)}\nonumber\\
&&\qquad \qquad \leqslant (1 + \alpha B)\,\spertxa{|g(X_T)|^2} + \left(\alpha B + \frac{1}{\beta} \right)\,\spertxa{\int_s^T |f(r,X_r,I_r)|^2 \,dr} \nonumber\\
&&\qquad \qquad + \frac{C^2}{\alpha}  + (\gamma + \beta)\, \spertxa{\int_s^T|Y_r^n|^2\,dr},
\quad s \in [t,T].\nonumber
\end{eqnarray}
Hence, by choosing $\alpha \in \left(0,\,\frac{1}{B}\right)$, $\gamma > \frac{T\,\lambda_0(A)}{1-\alpha B}$, $\beta >0$, and applying Gromwall's lemma to $s \rightarrow \spertxa{|Y_s^n|^2}$, we obtain:
\begin{eqnarray}\label{estim_Zn_Kn}
&& \sup_{ s \in [t,\, T]}\spertxa{|Y_{s}|^2} + \spertxa{\int_{t}^{T}\int_{E \times A}|Z^n_s(y,\,b)|^2 \, \tilde{p}(ds\,dy\, db)}\nonumber\\
&&\qquad \qquad \leqslant  C' \,\left(\spertxa{|g(X_T)|^2} + \,\spertxa{\int_t^T |f(s,X_s,I_s)|^2 \,ds} + C^2\,\right),
\end{eqnarray}
for some $C' >0$  depending only on $T$, which gives the required uniform estimate for $(Z^n)$ and also $(K^n)$ by \eqref{partial_estim_K^n}.
\endproof

We can finally present the
conclusion of the proof of Theorem \ref{Thm_ex_uniq_minimal_BSDE}:

\proof
Let $(t,x,a)\in[0,T]\times E\times A$. We first show that $(Y^n,Z^n,K^n)$ (we omit the dependence on $(t,x,a)$ for simplicity of notation) solution to \eqref{BSDE_penalized_K} converges  in a suitable way to some process $(Y,Z,K)$ solution to the  constrained BSDE \eqref{BSDE}-\eqref{BSDE_constraint}.
By Lemma \ref{Lemma_monotonicity_boundedness_Y_n}, $(Y^n)_n$ converges increasingly to some adapted process $Y$, which  moreover satisfies $\spertxa{\sup_{s \in [t,T]}\, |Y_s|^2}< \infty$
by the uniform estimate for $(Y^n)_n$ in Lemma \ref{lemma_BSDE_estimations} and Fatou's lemma.
Furthermore, by the dominated convergence theorem,
we also have $\E\int_0^T|Y^n_t-Y_t|^2dt\to 0$.
Next, we prove that there exists $(Z,\,K) \in \mathbb{\textbf{L}}^{\textbf{2}}_{\textbf{t,x,a}}\textup{(q)}\times \mathbb{\textbf{K}}^{\textbf{2}}_{\textbf{t,x,a}}$ with $K$ predictable, such that
\begin{itemize}
\item[(i)] $Z$ is the weak limit of $(Z^n)_n$ in $\mathbb{\textbf{L}}^{\textbf{2}}_{\textbf{t,x,a}}\textup{(q)}$;
\item[(ii)] $K_{\tau}$ is the weak limit of $(K_{\tau}^n)_n$ in $\mathbb{\textbf{L}}^{\textbf{2}}_{\textbf{t,x,a}}(\mathcal{F}_{\tau})$, for any stopping time $\tau$ valued in $[t\,\,T]$;
\item[(iii)] $\P^{t,x,a}$-a.s.,
\begin{eqnarray*}
Y_{s} &=& g(X_T) + \int_{s}^{T}f(r,X_r,I_r)\,dr + K_T - K_s \nonumber \\
&& - \int_{s}^{T}\int_{E \times A}Z_r(y,\,b)\, q(dr\,dy\, db)  - \int_{s}^{T}\int_{A}Z_r(X_r,\,b)\, \lambda_0(db)\,dr, \,\, s\in[t,T],
\end{eqnarray*}
with
\begin{equation*}
Z_s(X_{s-},b)\leqslant 0, \qquad  ds \otimes d\P^{t,x,a} \otimes \lambda_0(db)-\text{a.e.}\,\,%\text{on}\,\, [t,T]\times \Omega \times A.
\end{equation*}
\end{itemize}
Let define the following mappings from $\mathbb{\textbf{L}}^{\textbf{2}}_{\textbf{t,x,a}}\textup{(q)}$ to $\mathbb{\textbf{L}}^{\textbf{2}}_{\textbf{t,x,a}}(\mathcal{F}_{\tau})$:
\begin{eqnarray*}
&&I_{\tau}^1  : \quad Z \mapsto \int_t^\tau \int_{E\times A}\, Z_s(y,b)\, q(ds\,dy \,db),\\
&&I_{\tau}^2  : \quad Z \mapsto \int_t^\tau \int_{A }\, Z_s(X_s,b)\, \lambda_0(db)\,ds,
\end{eqnarray*}
for each $\mathbb{F}$-stopping time  $\tau$ with values in $[t,T]$.
We wish to prove  that $I_{\tau}^1 Z^n$ and  $I_{\tau}^2 Z^n$ converge weakly
in $\mathbb{\textbf{L}}^{\textbf{2}}_{\textbf{t,x,a}}(\mathcal{F}_{\tau})$
 to  $I_{\tau}^1 Z$ and $I_{\tau}^2 Z$ respectively. Indeed, by the uniform estimates for $(Z^n)_n$ in Lemma \ref{lemma_BSDE_estimations}, there exists a subsequence,  denoted $(Z^{n_k})_{k}$, which converges weakly in $\mathbb{\textbf{L}}^{\textbf{2}}_{\textbf{t,x,a}}\textup{(q)}$.
Since  $I_1$ and $I_2$ are linear continuous operators they are also weakly continuous
so that we have  $I_{\tau}^1 Z^{n_k}\to I_{\tau}^1Z$ and
$I_{\tau}^2 Z^{n_k}\to I_{\tau}^2 Z$
 weakly
in $\mathbb{\textbf{L}}^{\textbf{2}}_{\textbf{t,x,a}}(\mathcal{F}_{\tau})$
as $k\to\infty$.
 Since we have from \eqref{BSDE_penalized_K}
\begin{eqnarray*}
K_\tau^{n_k} &=&  -Y_{\tau}^{n_k} + Y_{t}^{n_k} - \int_{t}^{\tau} f(r,X_r,I_r)\, dr  \nonumber\\
&& + \int_{t}^{\tau} \int_{A} Z_{r}^{n_k}(X_{r},b) \,\lambda_{0}(db)\, dr + \int_{t}^{\tau} \int_{E \times A} Z_{r}^{n_k}(y,b) \, q(dr\,dy\,db),
\end{eqnarray*}
we also obtain the weak convergence in $\mathbb{\textbf{L}}^{\textbf{2}}_{\textbf{t,x,a}}(\mathcal{F}_{\tau})$
as $k\to\infty$
\begin{eqnarray}\label{K_intermediate_def}
K_\tau^{n_k} \rightharpoonup K_\tau &:=&  -Y_{\tau} + Y_{t} - \int_{t}^{\tau} f(r,X_r,I_r)\, dr  \nonumber\\
&& + \int_{t}^{\tau} \int_{A} Z_{r}(X_{r},b) \,\lambda_{0}(db)\, dr + \int_{t}^{\tau} \int_{E \times A} Z_{r}(y,b) \, q(dr\,dy\,db).
\end{eqnarray}
Arguing as in  \cite{Pe} proof of Theorem 2.1,
or  \cite{KhMaPhZh} Lemma 3.5,
\cite{Es}
Theorem 3.1
we see that
  $K$ inherits from $K^{n_k}$ the properties of having nondecreasing paths
  and of being square integrable and  predictable.
  Finally, from Lemma 2.2 in \cite{Pe} it follows that $K$ and $Y$ are càdlàg, so that $K^{t,x,a}\in \mathbf{K_{t,x,a}^2}$ and $Y^{t,x,a}\in\mathbf{S_{t,x,a}^2}$.

Notice that the processes $Z$ and $K$ in \eqref{K_intermediate_def} are uniquely determined. Indeed, if $(Z,K)$ and $(Z',K')$ satisfy  \eqref{K_intermediate_def}, then the predictable processes $Z$ and $Z'$ coincide at the jump times and can be identified almost surely
with respect to
$\tilde{p}(\omega,ds\,dy\,db)\P^{t,x,a}(d\omega)$
 (a similar argument can be found in the proof of Proposition \ref{Prop_uniq_min_sol} to which we refer for more details). Finally, recalling that the jumps of $p$ are totally inaccessible, we also obtain the uniqueness of the component $K$. The uniqueness of $Z$ and $K$ entails that all the sequences $(Z^n)_n$ and $(K^n)_n$ respectively converge (in the sense of points (i) and (ii) above) to $Z$ and $K$.

It remains to show that the jump constraint \eqref{BSDE_constraint} is satisfied. To this end, we consider the functional on $\mathbb{\textbf{L}}^{\textbf{2}}_{\textbf{t,x,a}}\textup{(q)}$ given by
\begin{eqnarray*}
G : \quad Z &\mapsto& \spertxa{\int_{t}^{T}\int_A \, [Z_s(X_{s-},b)]^+\, \lambda_0(db)\, ds}.
\end{eqnarray*}
From uniform estimate \eqref{unif_pen_estimate}, we see that
$G(Z^n)\rightarrow 0$ as $n \rightarrow \infty$. Since $G$ is convex and strongly continuous in the strong topology of $\mathbb{\textbf{L}}^{\textbf{2}}_{\textbf{t,x,a}}\textup{(q)}$, then $G$ is lower semicontinuous in the weak topology of $\mathbb{\textbf{L}}^{\textbf{2}}_{\textbf{t,x,a}}\textup{(q)}$, see, e.g., Corollary 3.9 in \cite{Bre}. Therefore, we find
\begin{displaymath}
G(Z)\leqslant \liminf_{n \rightarrow \infty}G(Z^n)=0,
\end{displaymath}
from which follows the validity of the jump constraint \eqref{BSDE_constraint} on $[t,T]$. We have then  showed that $(Y,Z,K)$ is a solution to the constrained BSDE \eqref{BSDE}-\eqref{BSDE_constraint}. It remains to prove that this is the minimal solution. To this end, fix $n \in \N$ and consider  a triple $(\bar{Y},\bar{Z},\bar{K}) \in \mathbb{\textbf{S}}^{\textbf{2}}_{\textbf{t,x,a}} \times \mathbb{\textbf{L}}^{\textbf{2}}_{\textbf{t,x,a}}\textup{(q)} \times \mathbb{\textbf{K}}^{\textbf{2}}_{\textbf{t,x,a}}$ satisfying \eqref{BSDE}-\eqref{BSDE_constraint}.
For any $\nu \in \mathcal{V}^n$, by introducing the compensated martingale measure $q^\nu$,
%(dt \, dy \,da)= q(dt \, dy \,da)- (\nu_t(a)-1)d_1(y,a)\tilde{p}(X_{s-},I_{s-}, dy \, da)\, ds$,
we see that the solution $(\bar{Y},\bar{Z},\bar{K})$ satisfies: $\P^{t,x,a}$-a.s.,
\begin{eqnarray}\label{BSDE_nu_1}
\bar{Y}_s &=& g(X_T) + \int_{s}^{T}f(r,X_r,I_r)\,dr + \bar{K}_T - \bar{K}_s  \\
&& - \int_{s}^{T}\int_{E \times A}\bar{Z}_r(y,\,b)\, q^{\nu}(dr\,dy\, db)  - \int_{s}^{T}\int_{A}\nu_r(b)\,\bar{Z}_r(X_r,\,b)\, \lambda_0(db)\,dr \quad s\in[t,T].\nonumber
\end{eqnarray}
By taking the expectation under $\P^{t,x,a}_\nu$ in \eqref{BSDE_nu_1}, recalling  Lemma \ref{lemma_P_nu_martingale}, and that $\bar{K}$ is nondecreasing, we have
\begin{eqnarray}
\label{dual_repr}
\bar{Y}_{s} &\geqslant& \spernu{g(X_T) + \int_{s}^{T}f(r,X_r,I_r)\,dr} - \spernu{\int_{s}^{T}\int_{A}\nu_r(b)\,\bar{Z}_r(X_r,\,b)\, \lambda_0(db)\,dr} \notag \\
&\geqslant& \spernu{g(X_T) + \int_{s}^{T}f(r,X_r,I_r)\,dr}\qquad s\in [t,T],
\end{eqnarray}
since $\nu$ is valued in $(0,\,n]$ and $Z$ satisfies constraint \eqref{BSDE_constraint}. As $\nu$ is arbitrary in  $\mathcal{V}^n$, we get from the representation formula \eqref{rep_Y_n} that $\bar{Y}_s \geqslant Y^n_s$, $\forall \, s \in [t,T]$, $\forall \, n \in \N$. In particular, $Y_s = \lim_{n \rightarrow \infty}Y_s^n \leqslant \bar{Y}_s $, i.e., the minimality property holds.
The uniqueness of the minimal solution straightly follows from Proposition \ref{Prop_uniq_min_sol}.

To conclude the proof, we argue on the limiting behavior of the dual representation for $Y^n$ when $n$ goes to infinity.
Since $\mathcal{V}^n \subset \mathcal{V}$, it is clear from the representation \eqref{rep_Y_n} that, for all $n$ and $s \in [t,T]$,
$Y_s^n \leqslant
\esssup_{\nu \in \mathcal{V}}\spernu{g(X_{T}) + \int_{s}^{T} f(r,X_{r},I_{r})\, dr
\,\bigg\vert\,\calf_s}$.
Moreover, being $Y$ the pointwise limit of $Y^n$, we deduce that
\begin{equation}\label{Y_ineq}
Y_s = \lim_{n \rightarrow \infty}Y_s^n \leqslant \esssup_{\nu \in \mathcal{V}}\spernu{g(X_{T}) + \int_{s}^{T} f(r,X_{r},I_{r})\, dr\,\bigg\vert\,\calf_s}.
\end{equation}
On the other hand, for any $\nu \in \mathcal{V}$, introducing  the compensated martingale measure $q^\nu$ under $\P^\nu$ as usual,
we see that $(Y,\,Z,\,K)$ satisfies
\begin{eqnarray}\label{BSDE_nu}
Y_{s} &=& g(X_T) + \int_{s}^{T}f(r,X_r,I_r)\,dr + K_T - K_s \\
&& - \int_{s}^{T}\int_{E \times A}Z_r(y,\,b)\, q^{\nu}(dr\,dy\, db)  - \int_{s}^{T}\int_{A} Z_r(X_r,\,b)\,\nu_r(b) \lambda_0(db)\,dr, \qquad s \in [t,T]. \notag
\end{eqnarray}
Arguing in the same way as in \eqref{dual_repr},
we obtain
\begin{eqnarray*}
Y_{s}
&\geqslant& \spernu{g(X_T) + \int_{s}^{T}f(r,X_r,I_r)\,dr\,\bigg\vert\,\calf_s},
\end{eqnarray*}
so that
$Y_s \geqslant \esssup_{\nu \in \mathcal{V}}\spernu{g(X_{T}) + \int_{s}^{T} f(r,X_{r},I_{r})\, dr
\,\Big\vert\, \calf_s}$
  by the arbitrariness of $\nu \in \mathcal{V}$.
  Together with \eqref{Y_ineq} this gives the required equality.
\endproof

\section{A BSDE representation   for the value function}
\label{FeynmanKac}

In this section we conclude the last step in the method of control randomization
and we show that the minimal solution to the constrained  BSDE
 \eqref{BSDE}-\eqref{BSDE_constraint} actually provides  a
 non-linear Feynman-Kac representation
 of the solution to the Hamilton-Jacobi-Bellman (HJB) equation \eqref{HJB}-\eqref{HJB_term},
 that we re-write here:
$$
-\frac{\partial v}{\partial t}(t,x)=\sup_{a \in A}\left(\mathcal{L}^{a}_E v(t,x) + f(t,x,a)\right) ,
\qquad
v(T,x)= g(x).
$$
As a consequence of the dual representation in
 Theorem \ref{Thm_ex_uniq_minimal_BSDE}
 it follows
 that the value function of the original optimal control  problem can be identified with the dual one, which in particular  turns out to be indepedent on the variable $a$.

For our result we need the following conditions:
\begin{eqnarray}
 \label{lambda_unif_finite_bis}
&& \sup_{x \in E, a \in A}\, \lambda(x,a,E) < \infty,
\\\label{Fellerbis}
&& \lambda  \text{ is a Feller transition kernel,}
\\
&& \label{fgcont}
f\in C_b([0,T]\times E\times A), \quad g\in C_b(E).
\end{eqnarray}
We note that these assumptions are stronger that those required in
Theorem \ref{thm_HJB_ben_posta}
and therefore they imply that there exists a unique solution
$v\in LSC_b([0,T]\times E)$
to
the HJB equation in the sense of
Definition \ref{def_sol_HJB}. If, in addition, $A$ is a compact
metric space then
$v\in C_b([0,T]\times E)$ by Corollary \ref{thm_HJB_ben_postaC}.

Let us consider again the Markov process $(X,I)$ in $E\times A$ constructed in
section \ref{X_I_construction}, with corresponding family of probability
measures $\P^{t,x,a}$
and generator
$\mathcal{L}$ introduced in
\eqref{generatoredoppio}.
Since
\eqref{lambda_unif_finite_bis}-\eqref{fgcont} are also stronger than
\textbf{\textup{(H$\lambda$)}} and \textbf{\textup{(Hfg)}}, by
Theorem \ref{Thm_ex_uniq_minimal_BSDE}   there exists a unique solution
to the BSDE
\eqref{BSDE}-\eqref{BSDE_constraint}.

Our main result is as follows:

\begin{theorem}\label{thm_conv_vn}
Assume
 \eqref{lambda_unif_finite_bis}, \eqref{Fellerbis},
 \eqref{fgcont}.
Let  $v$ be the unique solution to the Hamilton-Jacobi-Bellman equation
provided by Theorem \ref{thm_HJB_ben_posta}.
Then for every $(t,x,a) \in [0,\,T] \times E \times A$,
\begin{displaymath}
v(t,x) = Y_t^{t,x,a},
\end{displaymath}
 where $Y^{t,x,a}$ is the first component of the
  minimal solution to the constrained
  BSDE with non\-po\-si\-ti\-ve jumps \eqref{BSDE}-\eqref{BSDE_constraint}.

More generally, we have  $\P^{t,x,a}$-a.s.,
\begin{displaymath}
v(s,X_s) = Y_s^{t,x,a}, \qquad s \in [t,T].
\end{displaymath}
Finally, for the value function $V$ of the optimal control problem
 defined in  \eqref{value_function} and
 the dual value function $v^\ast$
 defined in \eqref{funzione_valore_duale}
 we have the equalities
\begin{displaymath}
V(t,x)=v(t,x) =Y_t^{t,x,a}= v^{\ast}(t,x,a).
\end{displaymath}
In particular, the latter functions do not depend on $a$.
\end{theorem}

The rest of this section is devoted to prove Theorem \ref{thm_conv_vn}.

\subsection{A  penalized HJB equation}

Let us recall the penalized BSDE associated to \eqref{BSDE}-\eqref{BSDE_constraint}: $\P^{t,x,a}$-a.s.,
\begin{eqnarray}
Y_{s}^{n,t,x,a} &=& g(X_{T}) + \int_{s}^{T} f(r,X_{r},I_{r})\, ds   - \int_{s}^{T} \int_{E \times A} Z_{r}^{n,t,x,a}(y,b) \, q(dr\,dy\,db)\label{BSDE_pen_rap}\\
&& + \int_{s}^{T} \int_{A}\left\{ n\,[Z_{r}^{n,t,x,a}(X_{r},b)]^{+} - Z_{r}^{n,t,x,a}(X_{r},b)\right\} \,\lambda_{0}(db)\, dr, \qquad s \in [t,T].\nonumber
\end{eqnarray}
Let us now consider the parabolic semi-linear penalized integro-differential
equation, of HJB type: for any $n \geq1$,
\begin{eqnarray}
\frac{\partial v^n}{\partial t}(t,x,a)+\mathcal{L} v^n(t,x,a) + f(t,x,a) \label{HJB_pen}\\
+\int_{A}\{n\,[v^n(t,x,b)-v^n(t,x,a)]^+\,-(v^n(t,x,b)-v^n(t,x,a))\}\,\lambda_0(db)&=& 0 \quad \text{on}\,\,[0,\,T) \times E\times A,\nonumber\\
v^n(T,x,a)&=& g(x) \quad \text{on}\,\,E\times A,\label{HJB_pen_term}
\end{eqnarray}
The following lemma states that
the solution of \eqref{HJB_pen}-\eqref{HJB_pen_term}
can be represented probabilistically by means of the solution to the penalized BSDE \eqref{BSDE_pen_rap}:
\begin{lemma}\label{thm_well_pos_vn}
Assume
 \eqref{lambda_unif_finite_bis}, \eqref{Fellerbis},
 \eqref{fgcont}.
 Then there exists a unique function
$ v^n\in C_b([0,T]\times E\times A)$ such that $t\mapsto v^n(t,x,a)$ is continuously
differentiable on $[0,T]$
and \eqref{HJB_pen}-\eqref{HJB_pen_term}
hold for every $(t,x,a)\in [0,\,T) \times E \times A$.

Moreover,  for every $(t,x,a)\in [0,\,T]\times E \times A$ and  for every $n \in \N$,
\begin{eqnarray}
Y_s^{n,t,x,a}&=&v^n(s,X_s,I_s)\label{Yn_id}\\
Z_s^{n,t,x,a}(y,b)&=&v^n(s,y,b)-v^n(s,X_{s-},I_{s-}),\label{Zn_id}
\end{eqnarray}
(to be understood as an equality between
elements
of the space  \emph{$\mathbb{\textbf{S}}^{\textbf{2}}_{\textbf{t,x,a}} \times \mathbb{\textbf{L}}^{\textbf{2}}_{\textbf{t,x,a}}\textup{(q)}$})
so that in particular $v^n(t,x,a)= Y^{n,t,x,a}_t$.
\end{lemma}
\proof
We first note that $ v^n\in C_b([0,T]\times E\times A)$ is the required solution
 if and only if
\begin{eqnarray}\label{HJB_pen_integral}
v^n(t,x,a) = g(x) +\int_{t}^{T}\,\mathcal{L} v^n(s,x,a)\,ds + \int_{t}^{T}\,f^n(s,x,a, v^n(s,x,\cdot)-v^n(s,x,a))
\end{eqnarray}
for $ t \in [0,,T)$, $ x\in E$, $a \in  A$,
where $f^n(t,x,a,\psi)$ is the map defined in \eqref{f^n}.
We  use a fixed point argument, introducing a map  $\Gamma$
from $C_b([0,T]\times E\times A)$
to itself setting $v=\Gamma(w)$ where
\begin{displaymath}
v(t,x,a) = g(x) + \int_{t}^{T}\mathcal{L} w(s,x,a)\,ds + \int_{t}^{T} f^n(s,x,a,w(s,x,\cdot)-w(s,x,a))\,ds.
\end{displaymath}
Using the boundedness assumptions on $\lambda$ and $\lambda_0$
it can be shown by standard arguments that some iteration of
the above map is a contraction in the space of bounded measurable real functions on $[0,\,T]\times E \times A$ endowed with the supremum norm and therefore
 the map $\Gamma$ has a unique fixed point, which is the required solution $v^n$.

We finally prove the identifications \eqref{Yn_id}-\eqref{Zn_id}.
Since $ v^n\in C_b([0,T]\times E\times A)$
 we can apply  the It$\hat{\mbox{o}}$ formula
 to the process $v(s,X_{s},I_{s})$, $s \in [t,T]$, obtaining, $\P^{t,x,a}$-a.s.,
\begin{eqnarray*}
v^n(s,X_s,I_s) &=& v^n(t,x,a) + \int_{t}^{s}\left( \frac{\partial v^n}{\partial r}(r,X_r,I_r)+ \mathcal{L}^I_r v^n(r,X_r,I_r)\right)\,dr\nonumber\\
&&+ \int_{t}^{s}\int_{E \times A}\left( v^n(r,y,b)-v^n(r,X_{r-},I_{r-})\right)\,q(dr\,dy\,db),\qquad s \in [t,T].
\end{eqnarray*}
Taking into account that $v^n$ satisfies \eqref{HJB_pen}-\eqref{HJB_pen_term} and that $(X,I)$ has piecewise constant trajectories, we obtain $\P^{t,x,a}$-a.s.,
\begin{displaymath}
\frac{\partial v^n}{\partial r}(r,X_r,I_r)+\mathcal{L} v^n(r,X_r,I_r) + f^n(r,X_r,I_r, v^n(r,X_r,\cdot)-v^n(r,X_r,I_r)) = 0,
\end{displaymath}
for almost all $r \in [t,T]$. It follows that, $\P^{t,x,a}$-a.s.,
\begin{eqnarray*}
v^n(s,X_s,I_s) &=& v^n(t,x,a) - \int_{t}^{s}f^n(r,X_r,I_r, v^n(r,X_r,\cdot)-v^n(r,X_r,I_r))\,dr\nonumber\\
&&+ \int_{t}^{s}\int_{E \times A}\left( v^n(r,y,b)-v(r,X_{r-},I_{r-})\right)\,q(dr\,dy\,db),\qquad s \in [t,T].
\end{eqnarray*}
Since $v^n(T,x,a) = g(x)$ for all $(x, a) \in E \times A$, simple passages show that
\begin{eqnarray*}
v^n(s,X_s,I_s) &=& g(X_T) + \int_{t}^{s}f^n(r,X_r,I_r, v^n(r,X_r,\cdot)-v^n(r,X_r,I_r))\,dr\nonumber\\
&&- \int_{t}^{s}\int_{E \times A}\left( v^n(r,y,b)-v(r,X_{r-},I_{r-})\right)\,q(dr\,dy\,db),\qquad s \in [t,T].
\end{eqnarray*}
Therefore the pairs $(Y_s^{n,t,x,a},Z_s^{n,t,x,a}(y,b))$ and $(v^n(s,X_s,I_s),v^n(s,y,b)-v^n(s,X_{s-},I_{s-}))$ are both solutions to the same BSDE under $\P^{t,x,a}$, and thus they  coincide as members of the space  $\mathbb{\textbf{S}}^{\textbf{2}}_{\textbf{t,x,a}} \times \mathbb{\textbf{L}}^{\textbf{2}}_{\textbf{t,x,a}}\textup{(q)}$. The required equalities \eqref{Yn_id}-\eqref{Zn_id} follow. In particular we have that $v^n(t,x,a) = Y^{n,t,x,a}_t$.
\endproof

\subsection{Convergence of the penalized solutions and conclusion of the proof}

We study the behavior of the functions $v^n$ as $n\to\infty$. To this end
we first show that they are  bounded above by the solution
to the HJB equation.

\begin{lemma}\label{lem_comp_sup_sol}
Assume
 \eqref{lambda_unif_finite_bis}, \eqref{Fellerbis},
 \eqref{fgcont}.
Let $v$ denote the solution to the HJB equation
 as provided by
Theorem \ref{thm_HJB_ben_posta}
and $v^n $  the solution to \eqref{HJB_pen}-\eqref{HJB_pen_term}
as provided in Lemma \ref{thm_well_pos_vn}. Then,
for all $(t,x,a)\in [0,\,T]\times E \times A$ and $n\ge 1$,
\begin{displaymath}
v(t,x) \geq  v^n(t,x,a).
\end{displaymath}
\end{lemma}

\proof
Let $v : [0,\,T] \times E  \rightarrow \R$ be a  solution to the
HJB equation.
As in the proof of
Proposition \ref{prop_construction_X_Markov_pure_jump_1} we have the identity
$$
g(X_T)- v(t,X_t) =   \int_{t}^T  \frac{\partial v}{\partial r}(r,X_r)\,dr
+  \int_{(t,T]}\int_{E\times A} (v(r,y) - v(r,X_{r-}))\,p(dr\,dy\,db),
$$
which follows from the absolute continuity of $t\mapsto v(t,x)$,
taking into account that $X$ is constant among
jump times and using the definition of the random measure $p$
defined in \eqref{misura_aleatoria_X_I}
and the fact that $v$ depends on $t,x$ only.
Since $v$ is a solution to the HJB equation we have,
for all $x  \in   E$ $a\in A$,
$$
- \frac{\partial v}{\partial t}(t,x) \geq
  \mathcal{L}^{a}_E v(t,x) + f(t,x,a)
   =\int_{E } (v(t,y) - v(t,x))\,\lambda(x,a,dy)
   + f(t,x,a),
$$
almost surely on $ [0,T]$.
Taking into account that  $(X,I)$ has piecewise constant trajectories we obtain
\begin{eqnarray}\label{v_soprasol_di_v_n}
  g(X_T)- v(t,X_t) &\le & \int_{(t,T]}\int_{E\times A} (v(r,y) - v(r,X_{r-}))\,p(dr\,dy\,db)
   \\\nonumber
    & & -\int_{t}^{T}  \int_{E } (v(r,y) - v(r,X_{r}))\,\lambda(X_r,I_r,dy)\,dr
    -\int_{t}^{T}   f(r,X_r,I_r)\,dr.
\end{eqnarray}
Then, for any $n \ge 1$  and $\nu \in \mathcal{V}^n$ let us
consider the probability $\P^{t,x,a}_\nu$ introduced above
 and recall that under $\P^{t,x,a}_\nu$ the  compensator of the random
 measure $p(dr\,dy\,db)$ is
$ \tilde{p}^\nu(dr\,dy\,db)=  \nu_r(b) \, \lambda_0 (db) \, \delta_{\{ X_{r-} \}}(dy) \,  dr+
\lambda(X_{r-},\,I_{r-},\,dy)\, \delta_{\{ I_{r-} \}}(db)\, dr$. Noting
that $v(r,y) - v(r,X_{r-})$ is predictable, taking the expectation in
\eqref{v_soprasol_di_v_n} we obtain
$$
\E^{t,x,a}_\nu[g(X_T)]- v(t,x) \le
    -\E^{t,x,a}_\nu\int_{t}^{T}   f(r,X_r,I_r)\,dr.
$$
Since   $\nu \in \mathcal{V}^n$ was arbitrary, and recalling
\eqref{rep_Y_n},
we conclude that
$$
v(t,x)   \geq  \sup_{\nu \in \mathcal{V}^n}\spernu{g(X_T)+ \int_{t}^{T}\,f(r,X_r,I_r)\,dr}
 =  v^n(t,x,a).
$$
\endproof

From
Lemma \ref{thm_well_pos_vn} we know that $v^n(t,x,a)= Y^{n,t,x,a}_t$, and
from
Lemma \ref{Lemma_monotonicity_boundedness_Y_n} we know that
$v^n(t,x,a)$ is monotonically increasing and uniformly bounded. Therefore we can define
$$
\bar v(t,x,a):=\lim_{n\to\infty }v^n(t,x,a),
\qquad t\in [0,T],\, x\in E,\,a\in A.
$$
$\bar v$ is bounded,
and from Lemma \ref{lem_comp_sup_sol} we deduce that  $\bar v\le v$.
As an increasing limit of continuous functions,
$\bar v$ is lower semi-continuous. Further properties of $\bar v$
are proved in the following lemma. In particular,
\eqref{bar_v_soprasol} (or
\eqref{bar_v_soprasol_a_variabile}) means that $\bar v$
is a supersolution to the HJB equation.

\begin{lemma}\label{lem_indip_da_a}
Assume
 \eqref{lambda_unif_finite_bis}, \eqref{Fellerbis},
 \eqref{fgcont}  and let $\bar v$ be the increasing limit
 of $v^n$. Then $\bar v$ does not depend on $a$, i.e.
 $\bar v(t,x,a)=\bar v(t,x,b)$ for every
 $t\in [0,T]$, $ x\in E$ and $a,b\in A$. Moreover, setting
 $\bar v(t,x)=\bar v(t,x,a)$ we have
 \begin{eqnarray}\label{bar_v_soprasol}
    \bar v(t,x)-\bar v(t',x) &\ge&  \int_t^{t'}(\call_E^a\bar v(s,x) +f(s,x,a))\,ds,
 \quad
 0\le t\le t'\le T,\, x\in E,\,a\in A.
 \end{eqnarray}
 More generally, for arbitrary Borel-measurable $\alpha:[0,T]\to A$ we have
 \begin{eqnarray}\label{bar_v_soprasol_a_variabile}
 \bar v(t,x)-\bar v(t',x) \ge \int_t^{t'}(\call_E^{\alpha(s)}\bar v(s,x) +f(s,x,\alpha(s)))\,ds,
 \qquad
 0\le t\le t'\le T,\, x\in E,\,a\in A.
\end{eqnarray}

\end{lemma}

\proof
$v^n$ satisfies the integral equation
  \eqref{HJB_pen_integral}, namely
\begin{eqnarray*}
v^n(t,x,a) &=& g(x) +
\int_{t}^{T}\int_E\, (v^n(s,y,a)-v^n(s,x,a))\,\lambda(x,a,dy)\,ds
\\&&
+\int_t^{T}f(s,x,a)\,ds
+ n\int_{t}^{T}\int_A\,[v^n(s,x,b)- v^n(s,x,a)]^+\,\lambda_0(db)\,ds.
\end{eqnarray*}
Since $v^n$ is a bounded sequence in $C_b([0,T]\times E\times A)$ converging
pointwise to $\bar v$, setting $t=0$, dividing by $n$ and letting
$n\to\infty$ we obtain
\begin{equation}\label{bar_v_al_limite}
    \int_{0}^{T}\int_A\,[\bar v(s,x,b)-\bar  v(s,x,a)]^+\,\lambda_0(db)\,ds=0.
\end{equation}
Next we claim that $\bar v$ is right-continuous in $t$ on $[0,T)$,
for fixed $ x\in E$, $a\in A$.
To prove this we first note that, neglecting the term with the positive part $[\ldots]^+$
we have
\begin{eqnarray}
  v^n(t',x,a)-  v^n(t,x,a)& \le&
-\int_{t}^{t'}\int_E\, (v^n(s,y,a)-v^n(s,x,a))\,\lambda(x,a,dy)\,ds
-\int_t^{t'}f(s,x,a)\,ds\nonumber
\\&\le& C_0(t'-t),\label{slope}
\end{eqnarray}
for some constant $C_0>0$ and for all
$0\le t\le t'\le T$ and $n\ge 1$, where we have used again the fact that
$v^n$ is uniformly bounded. Now fix $t\in [0,T)$. Since,
as already noticed, $\bar v$ is lower semi-continuous we have
$\bar v(t,x,a)\le\liminf_{s\downarrow t}\bar v (s,x,a)$. The required right continuity
follows if we can show  that $\bar v(t,x,a)\ge\limsup_{s\downarrow t}\bar v (s,x,a)$.
Suppose not. Then there exists $s_k\downarrow t$ such that
$\bar v(s_k,x,a)$ tends to some limit $l>\bar v(t)$. It follows that
$  \bar v (s_k,x,a)-\bar   v(t,x,a)>C_0(s_k-t)$ for some $k$ sufficiently large,
and therefore also
$  v^n(s_k,x,a)-  v^n(t,x,a)>C_0(s_k-t)$ for some $n$ sufficiently large,
contradicting \eqref{slope}. This contradiction shows
that $\bar v$ is right-continuous in $t$ on $[0,T)$.

Then it follows from \eqref{bar_v_al_limite} that
$\int_A[\bar v(t,x,b)-\bar  v(t,x,a)]^+\,\lambda_0(db)=0$ for every
$x\in E$, $a\in A$, $t\in [0,T]$. Therefore there exists  $B\subset A$
(dependent on $t,x,a$) such that $B$ is a Borel set  with
$\lambda_0(B)=0$, and
\begin{eqnarray}
  \bar v(t,x,a) &\ge &  \bar v(t,x,b'),\qquad b'\notin B.\label{bar_v_ineq}
\end{eqnarray}
Since $\lambda_0$ has full support, $B$ cannot contain any open ball.
So given an arbitrary $b\in A$ we can find a sequence $b_n\to b$, $b_n\notin B$.
Writing \eqref{bar_v_ineq} with $b_n$  instead of $b'$ and using the
lower semi-continuity of $\bar v$ we deduce that
$\bar v(t,x,a) \ge  \liminf_n\bar v(t,x,b_n)\ge \bar v(t,x,b)$. Since $a$ and $b$ were
arbitrary we finally conclude that
 $\bar v(t,x,a)=\bar v(t,x,b)$ for every
 $t\in [0,T]$, $ x\in E$ and $a,b\in A$, so that $\bar v(t,x,a)$ does not depend
 on $a$ and we can define
 $\bar v(t,x)=\bar v(t,x,a)$.

Passing to the limit as $n\to\infty$ in the first inequality of \eqref{slope}
we immediately obtain
\eqref{bar_v_soprasol}, so it remains to prove \eqref{bar_v_soprasol_a_variabile}.
Let $\cala(\bar v)$ denote
the set of all   Borel-measurable $\alpha:[0,T]\to A$ such that
\eqref{bar_v_soprasol_a_variabile} holds, namely for every
$0\le t\le t'\le T$, $ x\in E$, $a\in A$,
\begin{eqnarray}\label{bar_v_soprasol_al_limite}
  \bar v(t,x)-\bar v(t',x) &\ge&
   \int_t^{t'}\int_E  \bar v(s,y)\,\lambda(x,\alpha(s),dy)\,ds
    \\&&
    -
\int_t^{t'}  \bar v(s,x)\,  \lambda(x,\alpha(s),E)\,ds
  +\int_t^{t'}f(s,x,\alpha(s))\,ds.
\end{eqnarray}
Suppose that $\alpha_n\in \cala(\bar v)$,
$\alpha:[0,T]\to A$
is Borel-measurable and $\alpha_n(t)\to \alpha_n(t)$
for almost all $t\in[0,T]$.  Note that
\begin{eqnarray}\label{lim_int_v_n}
\int_E  \bar v(t,y)\,\lambda(x,a,dy)=\lim_{n\to\infty}
\int_E  \bar v^n(t,y,a)\,\lambda(x,a,dy)
\end{eqnarray}
and the latter is an increasing limit. Since $v^n\in C_b([0,T]\times E\times A)$
and $\lambda$ is Feller,  for any $n\ge 1$
the functions in the right-hand side of \eqref{lim_int_v_n} are
continuous in $(t,x,a)$ (see e.g. \cite{BS}, Proposition 7.30) and therefore
the left-hand side is a lower semicontinuous function of $(t,x,a)$.
It follows from this and the Fatou Lemma that
\begin{eqnarray*}
\int_t^{t'}\int_E  \bar v(s,y)\,\lambda(x,\alpha(s),dy)\,ds
&\le&
\int_t^{t'}\liminf_{n\to\infty}\left[\int_E  \bar v(s,y)\,\lambda(x,\alpha_n(s),dy)\right]\,ds
\\&
\le &\liminf_{n\to\infty}
\int_t^{t'}\int_E  \bar v(s,y)\,\lambda(x,\alpha_n(s),dy)\,ds .
\end{eqnarray*}
Using this inequality and the continuity and boundedness of the maps
$a\mapsto \lambda(x,a,E)$,
 $a\mapsto  f(t,x,a)$ we see that
assuming the validity of
 inequality \eqref{bar_v_soprasol_al_limite} for $\alpha_n$ implies
 that it also holds for $\alpha$, hence $\alpha\in \cala(\bar v)$.

 Next we note that $\cala(\bar v)$ contains all piecewise  constant
 functions of the form $\alpha(t)$ $=$ $\sum_{i=1}^k $ $a_i 1_{[t_i,t_{i+1})}(t)$
 with $k\ge1$, $0=t_1<t_2<\ldots<t_{k+1}=T$, $a_i\in A$: indeed,
 it is enough to write down \eqref{bar_v_soprasol} with
$[t,t')=[t_i,t_{i+1})$ and sum up over $i$ to get
\eqref{bar_v_soprasol_a_variabile} for $\alpha(\cdot)$ and
therefore conclude that $\alpha(\cdot)\in\cala(\bar v)$. Since we have
already proved that the class $\cala(\bar v)$ is stable under
almost sure pointwise limits it follows that $\cala(\bar v)$
contains all
 Borel-measurable functions $\alpha:[0,T]\to A$ as required.
\endproof

We are now ready to conclude the proof of  our main result.

\noindent
{\it Proof of Theorem \ref{thm_conv_vn}.}
We will prove the inequality
\begin{eqnarray}
\label{v_bar_e_V}
  \bar v(t,x) &\ge & V(t,x),\qquad t\in [0,T] ,x\in E,
\end{eqnarray}
where $\bar v=\lim_{n\to\infty }v^n$ was introduced before Lemma
\ref{lem_indip_da_a}.
Since we know that $\bar v\le v$ and,
by
Theorem \ref{thm_verifica},  $v=V$ it follows from \eqref{v_bar_e_V} that
$\bar v=v=V$.
Passing to the limit as $n\to\infty$ in \eqref{Yn_id} and recalling
  \eqref{vstar_rep}
  all the other equalities follow immediately.

To prove \eqref{v_bar_e_V} we fix $t\in[0,T]$,  $x\in E$ and a
Borel-measurable map $\alpha: [0,T]\times E\to A$, i.e.
an element of $\cala_{ad}$, the set of
   admissible control laws for the primal control problem,
and denote by
$\mathbb{P}^{t,x}_{\alpha}$  the associated probability measure
on $(\Omega,\calf_\infty)$,
for the controlled system started at
time $t$ from point $x$, as in section \ref{Sect_control_problem}. We will prove that
$\bar v(t,x)\ge J(t,x,\alpha)$, the
gain functional defined in \eqref{functional_cost}.
Recall that in $\Omega$ we had defined a canonical marked point process
$(T_n,E_n)_{n\ge1}$ and the associated random measure $p$.
 Fix $\omega\in\Omega$ and consider the points
$T_n(\omega)$ lying in $(t,T]$, which we rename $S_i$;
thus, $t<S_1<\ldots S_k\le T$, for some $k$ (also depending on $\omega$).
Recalling that $\bar v(T,x)=g(x)$ we have
\begin{eqnarray*}
  g(X_T)-\bar v(t,x) &=& g(X_T)-\bar v(S_k,X_{S_k})
  +\sum_{i=1}^{k} [ \bar v(S_i,X_{S_i})- \bar v(S_i,X_{S_i-})]\\
  &&  +\sum_{i=2}^{k} [\bar v(S_i,X_{S_i-})-\bar v(S_{i-1},X_{S_{i-1}}) ]
  + \bar v(S_1,X_{S_1-})-\bar v(t,x).
\end{eqnarray*}
$\mathbb{P}^{t,x}_\alpha$-a.s we have  $X_{S_i-}= X_{S_{i-1}}$
($2\le i\le k$)
and $X_{S_1-}=x$, so we obtain
\begin{eqnarray*}
  g(X_T)-\bar v(t,x) &=& g(X_T)-\bar v(S_k,X_{S_k})
  +\sum_{i=1}^{k} [ \bar v(S_i,X_{S_i})- \bar v(S_i,X_{S_i-})]\\
  &&  +\sum_{i=2}^{k} [\bar v(S_i,X_{S_{i-1}})-\bar v(S_{i-1},X_{S_{i-1}}) ]
  + \bar v(S_1,x)-\bar v(t,x).
\end{eqnarray*}
The first sum can be written as
$$
\sum_{i=1}^{k} [ \bar v(S_i,X_{S_i})- \bar v(S_i,X_{S_i-})]=
\int_t^T\int_E[ \bar v(s,y)- \bar v(s,X_{s-})]\,p(ds\,dy),
$$
while the other can be estimated from above by repeated applications
of \eqref{bar_v_soprasol_a_variabile}, taking into account
that $X$ is constant in the intervals
$(t,S_1]$,  $(S_{i-1},S_{i}]$ ($2\le i\le k$) and
$(S_k,T]$:
\begin{eqnarray*}
\bar v(S_i,X_{S_{i-1}})-\bar v(S_{i-1},X_{S_{i-1}})
&\le&
-\int_{S_{i-1}}^{S_i}\Big(\call_E^{\alpha(s, X_{S_{i-1}})}
\bar v(s,X_{S_{i-1}}) +f(s,X_{S_{i-1}},\alpha(s,X_{S_{i-1}}))\Big)\,ds
\\
&=&
-\int_{S_{i-1}}^{S_i}\Big(\call_E^{\alpha(s, X_{s})}
\bar v(s,X_{s}) +f(s,X_{s},\alpha(s,X_{s}))\Big)\,ds
\end{eqnarray*}
for $2\le i\le k$
and similar formulae for the intervals  $(t,S_1]$, and
$(S_k,T]$. We end up with
\begin{eqnarray*}
 g(X_T)-\bar v(t,x)
&\le&
\int_t^T\int_E[ \bar v(s,y)- \bar v(s,X_{s-})]\,p(ds\,dy)
\\&&-
\int_{t}^{T}\Big(\call_E^{\alpha(s, X_{s})}
\bar v(s,X_{s}) +f(s,X_{s},\alpha(s,X_{s}))\Big)\,ds.
\end{eqnarray*}
Recalling that
the compensator of the measure $p$ under $\mathbb{P}^{t,x}_{\alpha}$ is
$1_{[t,\infty)}(s)\lambda(X_{s-},\,\alpha(s,X_{s-}),\,dy) \, ds$ we have, taking expectation,
$$
 \mathbb{E}^{t,x}_{\alpha}\int_t^T\int_E[ \bar v(s,y)- \bar v(s,X_{s-})]\,p(ds\,dy)
= \mathbb{E}^{t,x}_{\alpha}\int_{t}^{T}\call_E^{\alpha(s, X_{s})}
\bar v(s,X_{s})\,ds,
$$
which implies, by the previous inequality,
$\mathbb{E}^{t,x}_{\alpha}[
g(X_T)]-\bar v(t,x)\le -\mathbb{E}^{t,x}_{\alpha}
\int_{t}^{T}f(s,X_{s},\alpha(s,X_{s}))\,ds$ and so
$v(t,x)\ge
J(t,x,\alpha)$. Since $\alpha\in\cala_{ad}$ was
arbitrary we conclude that $v(t,x)\ge V(t,x) $.

\qed


\begin{thebibliography}{11}


%%\bibitem{AuFr}  Aubin, J.-P., Frankowska, H.
%% Set-valued analysis. Systems \& Control: Foundations \& Applications, 2. Birkh\"auser, 1990.
%

\bibitem{BaCo} Bandini, E., Confortola F.
Optimal control of semi-Markov processes with a backward stochastic differential equations approach. Preprint, arXiv: 1311.1063.


%\bibitem{BaBuPa} Barles, G., Buckdahn R., Pardoux E.
%Backward stochastic differential equations and integral-partial
%differential equations. Stochastics Stochastics Rep. 60 (1997),
%no. 1-2, 57-83.
%

\bibitem{Bech} Becherer, D.
Bounded solutions to backward SDE's with jumps for utility optimization and indifference hedging.
Ann. Appl. Probab.  16  (2006),  no. 4, 2027-2054.

%%\bibitem{Bis} Bismut, J.-M. A generalized formula of It\^{o} and some
%%other properties of stochastic flows. Z. Wahrsch. Verw. Gebiete 55
%%(1981), no. 3, 331-350.
%

\bibitem{BS} Bertsekas D.P., Shreve, S.E.
 Stochastic optimal control. The discrete time case. Mathematics in Science and Engineering, 139.
 Academic Press,  1978.


%\bibitem{Bo-Va-Wo} Boel, R., Varaiya, P., Wong, E.  Martingales on jump processes;
%Part I: Representation results; Part II: Applications; SIAM J. Control 13,
%999-1061.

\bibitem{BrLa}
Brandt, A., Last, G.  Marked point processes on the real line. The
dynamic approach. Springer, 1995.

%\bibitem{B} Br\'emaud,  P. Point processes and queues, Martingale dynamics.
%Springer Series in Statistics. Springer,  1981.

\bibitem{Bre} Brezis,  H. Functional Analysis. Sobolev Spaces and Partial Differential Equations. Springer,  2010.

  \bibitem{ChSoToVi07} Cheridito, P., Soner, M., Touzi, N.  Victoir, N.
    Second-order backward stochastic differential equations and fully nonlinear PDEs.
       Communication in Pure and Applied Mathematics
{\bf 60} (2007), 1081-1110.


%
%
%%\bibitem{CaFeSa} Carbone, R., Ferrario, B., Santacroce, M.
%%Backward stochastic differential equations driven
%%by c\`adl\`ag martingales. Theory   probab.   appl.   52 (2008), no.
%%2, 304-314.
%
%\bibitem{chito} Chitopekar, S. S.  Continuous time Markovian sequential control processes. SIAM J. Control, 7:367-389, 1969.
%
%\bibitem{Coh-Ell-1} Cohen, S. N., Elliott, R. J. Solutions of backward stochastic
%differential equations on Markov chains. Communications on Stochastic
%Analysis, 2(2):251-262, August 2008.
%
%\bibitem{Coh-Ell-2} Cohen, S. N., Elliott, R. J. Comparisons for backward stochastic
%differential equations on Markov chains and related no-arbitrage conditions.
%The Annals of Applied Probability, 20(1): 267-311, 2010.
%
%
%%\bibitem{Coh-Ell-3} Cohen, S.N. and Elliott, R.J. Existence, Uniqueness and Comparisons for BSDEs in General Spaces, to appear in Annals of Probability.
%

\bibitem{CoFu-mpp} Confortola, F., Fuhrman, M. Backward stochastic differential equations and optimal control of marked point processes. SIAM J. Control Optimization
     {\bf 51} (2013), no. 5, 3592-3623.


\bibitem{CoFu-m} Confortola, F., Fuhrman, M. Backward stochastic differential equations associated to jump Markov processes and their applications. Stochastic Processes and their Applications
    {\bf 124} (2014), pp. 289-316.

\bibitem{Co-Chou}  Choukroun, S., Cosso, A. Backward SDE Representation for Stochastic Control Problems with Non Dominated Controlled Intensity. Preprint arXiv 1405.3540  (2014).

\bibitem{CoFuPh} Cosso, A., Fuhrman, M., Pham, H.
Long time asymptotics for fully nonlinear Bellman equations:
a Backward SDE approach. Preprint arXiv 1410.1125  (2014).



%\bibitem{Cu-Fi-Ma-Te} Cuchiero, C., Filipovic, D., Mayerhofer, E., and Teichman, J. Affine Processes on Positive Semidefinite Matrices. Annals of Applied Probability 21 (2011), pp. 397-463.

%\bibitem{Da} Davis, M.H.A. Markov models and optimization.
%Monographs on Statistics and Applied Probability 49 (1993), Chapman and Hall, London.

%%\bibitem{Da-art} Davis, M.H.A. The representation of martingales of jump processes.
%%SIAM J. Control Optimization 14 (1976), no. 4, 623-638.
%
%\bibitem{Da-bo} Davis, M.H.A. Markov models and optimization.
%Monographs on Statistics and Applied Probability, 49, Chapman $\&$
%Hall,  1993.
%
%\bibitem{Da-Fa} Davis, M.H.A., Farid, M. Piecewise deterministic processes and viscosity solutions.
%McEneaney, W. M. et al. (ed) Stochastic Analysis, Control Optimization and Applications. A Volume in Honour of W. H. Fleming on Occasion of His 70th Birthday, Birkh\"auser (1999), 249-268.
%
%\bibitem{Dem} Dempster, M. A. H.
%Optimal control of piecewise deterministic Markov processes. Applied stochastic analysis (London, 1989), 303ï¿½325,
%Stochastics Monogr., 5, Gordon and Breach, New York, 1991.
%
%
\bibitem{ElKh10}
Elie, R., Kharroubi, I.
 Probabilistic representation and approximation for coupled systems of variational inequalities. Statist. Probab. Lett.  80  (2010),  no. 17-18, 1388-1396.


\bibitem{ElKh14} Elie, R., Kharroubi, I.
Adding constraints to BSDEs with jumps: an alternative to multidimensional reflections.
ESAIM Probab. Stat.  18  (2014), 233-250.

\bibitem{ElKh14a}
 Elie, R., Kharroubi, I.
 BSDE representations for optimal switching problems with controlled volatility. Stoch. Dyn.  14  (2014),  no. 3, 1450003, 15 pp.

\bibitem{EPQ}
El Karoui, N., Peng, S., Quenez, M.C. Backward stochastic
differential equations in finance. Math. Finance 7 (1997), no. 1,
1-71.

%\bibitem{ElK} El Karoui, N.
%Les aspects probabilistes du contr\^ole stocastique. [The probabilistic aspects of stochastic control] Ninth Saint Flour Probabilistic Summer School,-1979. (Saint Flour, 1979), pp. 73-238, Lecture Notes in Math., 876, Springer, 1981.
%
%% \bibitem{ElKHua} El Karoui, N., Huang, S.-J.
%%A general result of existence and uniqueness of backward
%%stochastic differential equations, in: Backward Stochastic
%%Differential Equations, N. El Karoui and L. Mazliak eds., Longman,
%%Harlow, 1997,   27ï¿½36.
%
%\bibitem{E} Elliott, R.J. Stochastic Calculus and its Applications. Springer,
%1982.
%
\bibitem{Es} Essaky, E. H. Reflected backward stochastic differential equation with jumps and RCLL obstacle. Bull. Sci. Math.  132  (2008),  no. 8, 690-710.


%\bibitem{G-S-II} Gihman I. I., Skorohod A. V.
%The Theory of Stochastic Processes II.
%Springer-Verlag Berlin Heidelberg New-York Tokyo, 1983.

%\bibitem{G-S-I} Gihman I. I., Skorohod A. V.
%The Theory of Stochastic Processes I.
%Springer-Verlag Berlin Heidelberg New-York Tokyo, 1983.

%\bibitem{Guo-He} Guo X., Hernandez-Lerma O.
%Continuous-time Markov decision processes: theory and applications.
%Stochastic Modelling and Applied Probability, 62. Springer-Verlag, Berlin, 2009.

%\bibitem{Howard} Howard, R. A. Dynamic Probabilistic Sistems.
%John Wiley, New York 1971.

%\bibitem{I-W} Ikeda N., Watanabe S.
%Stochastic differential equations and diffusion processes.
%North-Holland Mathematical Library, 24. North-Holland, 1989.

%\bibitem{J} Jacod, J. Calcul Stochastique et Problèmes de Martingales. Lecture Notes in Mathematics,
%Springer, Berlin  (1979).

\bibitem{FuPh} Fuhrman, M., Pham, H.
Randomized and Backward SDE representation for
optimal control of non-markovian SDEs. To appear on
The Annals of Applied Probability.

\bibitem{GuHe}
Guo, X.,   Hern\'andez-Lerma, O.
Continuous-time Markov decision processes.
Theory and applications. Stochastic Modelling and Applied Probability, 62. Springer, 2009.

\bibitem{J_a} Jacod, J. Multivariate point processes: predictable projection,
Radon-Nikodym derivatives, representation of martingales.
Z. Wahrscheinlichkeitstheorie und Verw. Gebiete 31 (1974/75), 235-253.

\bibitem{J-S} Jacod, J., Shiryaev, A.N. Limit theorems for stochastic processes. Second edition.
Springer,  1975.
%
%\bibitem{Jewell} Jewell, W. S. Markov-renewall programming I and II.
%Operational Res. 11 (1963), pp. 938-971, 938-971.
%


 \bibitem{KhMaPhZh}
 Kharroubi, I.,  Ma, J., Pham, H., Zhang J.
 Backward SDEs with contrained jumps and quasi-variational inequalities.
 Ann. Probab. 38
 (2010), no. 2, 794-840.

\bibitem{KhPh}
Kharroubi, I.,  Pham, H.
Feynman-Kac representation for Hamilton-Jacobi-Bellman IPDE.
To appear on  Ann. Probab.


\bibitem{KhLaPh}
 Kharroubi, I., Langren\'e, N., Pham, H.
A numerical algorithm for fully nonlinear HJB equations: an approach by control randomization. Monte Carlo Methods Appl.  20  (2014),  no. 2, 145-165.

\bibitem{KhLaPha}
Kharroubi, I., Langren\'e, N., Pham, H.
Discrete time approximation of fully nonlinear HJB equations via BSDEs with nonpositive jumps
Preprint  arXiv:1311.4505.


%%
%%\bibitem{K} Kunita, H. Stochastic differential equations and stochastic
%%flows of diffeomorphisms. \'Ecole d'\'et\'e de probabilit\'es de
%%Saint-Flour, XII-1982, 143-303, Lecture Notes in Math. 1097,
%%Springer,  1984.
%%
%
%
%
%\bibitem{Osaki} Osaki S., Mine H. Linear programming algorithms for semi-Markovian decision processes. J. Math. Anal. 22 (1968), pp. 356-381.
%
%
%\bibitem{PaPe} Pardoux E., Peng S. Adapted solution of a
%backward stochastic differential equation, Systems Control Lett.
%14 (1990) 55-61.
%%
%%\bibitem{Pe} Peng S. Stochastic Hamilton-Jacobi-Bellman equations.
%% SIAM J. Control Optim. 30 (1992),
%%284-304.

\bibitem{Pe} Peng, S. Monotonic limit theorem for BSDEs and non-linear Doob-Meyer decomposition.
Probability Theory and Related Fields. 16 (2000), pp. 225-234.


\bibitem{Pe06} Peng S.   $G$-expectation, $G$-Brownian motion
and related stochastic calculus of It\^o type.
 Stochastic analysis and applications,  541-567, Abel Symp., 2, Springer, Berlin, 2007.



\bibitem{Pl} Pliska, S.R. Controlled jump processes.
Stochastic Processes and their Applications 3 (1975), 259-282.

\bibitem{Pro-Shi} Protter, P., Shimbo, K. No arbitrage and general semimartingales. IMS Collections, Markov Processes and Related Topics: A Festschrift for Thomas G. Kurtz. 4 (2008), pp. 267-283.
%%
%\bibitem{Re-Yor} Revuz, D.,  Yor, M. Continuous Martingales and Brownian Motion.
%Grundlehren der mathematischen Wissenschaften. Springer, third
%edition, 1999.

%\bibitem{R-W-II} Rogers, L. C. G.,  Williams, D. Diffusions, Markov processes, and Martingales.
%Cambridge University Press, second edition, Volume 2, (2000).
%%
%
%
%\bibitem{Ross} Ross, S. M. Applied Probability Models with Optimization Applications. Holden-Day, San Francisco, 1970.
%
%\bibitem{Roy} Royer, M. Backward stochastic differential equations with jumps and
%related nonlinear expectations. Stochastic Processes and their Applications, 116(10):1358-1376, 2006.
%%
%%


\bibitem{SoToZh11} Soner M., Touzi N., and J. Zhang.
The wellposedness of second order backward SDEs.   Probability Theory and Related Fields
{\bf 153} (2011), 149-190.

% \bibitem{St1} Stone, L. D.  Necessary and sufficient conditions for optimal control of semi-Markov jump % processes. SIAM J. Control Optim. 11 (1973), no. 2, 187-201.
%
%%\bibitem{St2} Stone, L. D.  On the distribution of the supremum functional for the semi-Markov processes with continuous state spaces. The Annals of Mathematical Statistics 40 (1969), no. 3, 844-853.
%%
%\bibitem{TaLi}
%Tang S.,  Li X. Necessary Conditions for Optimal Control of
%Systems with Random Jumps. SIAM J. Control Optim. 32 (1994),
%1447-1475.
%
%
%\bibitem{Ver} V*ermes, D. Optimal control of piecewise deterministic Markov process. Stochastics 14 (1985), no. 3, 165ï¿½207.
%
%%
%%\bibitem{W} Wentzell, A.D. On the equation of the Theory of Conditional Markov Processes,
%% Theory  probab. appl.   10 (1965),  357-361.
%%
\bibitem{Xia} Xia, J. Backward stochastic differential equations with random measures.
Acta Mathematicae Applicatae Sinica 16 (2000), no. 3, 225-234.


\end{thebibliography}
\end{document}